\documentclass[a4paper,12pt, reqno]{amsart}
\usepackage[margin=2.5cm]{geometry}

\usepackage{graphicx}
\usepackage{amssymb,amsmath,amsthm}
\usepackage{float,geometry}

\usepackage{natbib}
\usepackage[hidelinks,colorlinks = true,linkcolor={blue!80!black},citecolor={blue!50!black},urlcolor={blue!80!black}]{hyperref}
\usepackage{cleveref} 

\usepackage{multirow}
\usepackage{multicol}
\usepackage{adjustbox}
\usepackage{booktabs}
\usepackage{subcaption}
\usepackage{pgfplots}

\usetikzlibrary{shapes,arrows}
\usepackage{setspace}
\usepackage{soul}

\newcommand{\dsum}{\displaystyle\sum}

\let\origmaketitle\maketitle
\def\maketitle{
	\begingroup
	\def\uppercasenonmath##1{} 
	\let\MakeUppercase\relax 
	\origmaketitle
	\endgroup
}

\begin{document}

\title[]{\Large A Novel Co-Evolutionary Algorithm for Solving a\\Bilevel Pricing and Hubs Location Problem\\under a Tree Topology}

\author[V. Blanco, J-F. Camacho-Vallejo, \MakeLowercase{and} C. Corpus]{
{\large V\'ictor Blanco$^{\dagger}$, Jos\'e-Fernando Camacho-Vallejo$^{\ddagger}$, and Carlos Corpus$^{\star}$}\medskip\\
$^\dagger$Institute of Mathematics (IMAG), Universidad de Granada, Spain\\
$^\ddagger$Tecnologico de Monterrey, Escuela de Ingenier\'ia y Ciencias, Monterrey, M\'exico\\
$^\star$Centro de Investigación en Ciencias Físico-Matemáticas (CICFIM), Universidad Autónoma de Nuevo León, M\'exico
\texttt{vblanco@ugr.es}, \texttt{fernando.camacho@tec.mx}, \texttt{carlos.ed240896@gmail.com}}

\maketitle

\begin{abstract}
This paper introduces the Bilevel Tree-of-Hubs Location Problem with Prices (BTHLPwP). The BTHLPwP is a multiple-allocation hub location problem in which, in addition to determining the nodes and links of a tree-shaped hub backbone network, the prices for using this network must also be set. We assume that two different types of agents make decisions in this problem. On the one hand, one agent (the leader) determines the structure and sets the prices for using the hub {backbone} network. On the other hand, the other agent (follower) decides on the optimal usage of the network. {The leader seeks to maximize its profit, while the follower aims to minimize the costs incurred for using the network to ship their commodities. We present a bilevel optimization formulation for this problem, followed by an equivalent single-level reformulation.} Then, we propose a novel Co-Evolutionary Algorithm (Co-EA) to solve three well-known datasets of instances adapted for our problem. The main novelty of the proposed Co-EA lies in the way the co-evolving populations are considered. While traditionally one population focuses on the leader's solutions and the other on the follower's, in our approach, each population is associated with a subset of the leader's decision variables. Consequently, the follower's optimal reaction is obtained for a specific decision made by the leader, resulting in bilevel feasible solutions. We then analyze the results obtained from extensive computational experimentation using the proposed Co-EA. {To validate its performance, we use the best solution found by the Co-EA as a warm-start solution in the single-level reformulation, which is then solved by a general-purpose solver. Additionally, we provide interesting managerial insights regarding the increase in the number of hubs in the backbone network.} This research concludes with recommendations for further research on this topic.
\end{abstract}


\section{Introduction} \label{sec:intro}

Globalization demands dynamic commercial connections that ensure the reliable delivery of commodities to users worldwide. The adequate design of commercial networks is crucial for both companies and users, enhancing their convenience. This type of decision typically involves addressing Hub Location Problems (HLPs), which integrate tools from Network Design and Location Science. In HLPs, the objective is to determine the optimal placement of certain special nodes in a network (hub nodes), which serve as transshipment points for routing commodities. Hub backbone networks provide practical benefits from an operational perspective, as they concentrate flow within a smaller network (the hub {backbone} network), facilitating efficient delivery monitoring and reducing users' routing costs. HLPs have garnered significant attention from researchers over the past decade due to their practical relevance and the mathematical and computational challenges they present (see \cite{campbell1994integer,contreras2019hub} for further details on hub location). {Recently, HLPs have been applied in various systems, including air transportation \citep{farahani2013hub}, postal services \citep{bagherinejad2020dynamic}, cargo delivery \citep{khaleghi2021robust}, and supply chains for perishable goods \citep{eydi2025hierarchical}.}  

As mentioned earlier, HLPs integrate decisions from Location Science, which determines the optimal positions of hub nodes and links, and Network Design, which facilitates the routing of commodities through the hub {backbone} network. Location Science focuses on identifying the optimal placement of services to adequately meet user demand. A substantial body of literature addresses various location problems and their solution methodologies. Prominent problems include the Weber Problem \cite{weber1929theory}, the Uncapacitated Facility Location Problem \cite{Cornue}, the Maximal Covering Location Problem \cite{church1974maximal}, and the $p$-median problem \cite{hakimi1983locating}. When the placement of services is influenced by routing costs associated with sending flows among users, location problems transition into HLPs, which have received extensive attention in recent years, emerging as an active research area \cite{HLP1,blanco2022hub,blanco2023hub,hekmatfar2009hub,nickel2001hub}. Most HLPs assume a complete hub {backbone} network, though this assumption is often unrealistic when the costs of installing links between hub nodes are disproportionately high. This discrepancy has spurred the exploration of alternative hub backbone network topologies, such as star-star \cite{labbe2008solving}, cycle \cite{contreras2017exact}, tree \cite{contreras2009tight,contreras2010tree}, or general shapes \cite{o2015hub,o2015multiple}. In the Tree-of-Hubs Location Problem (THLP), the hub backbone network is assumed to have a tree structure. This topology is particularly relevant in HLP scenarios where the setup costs for hub links exceed the routing costs, making it advantageous to construct a minimally connected network linking the activated hub nodes.

In this paper we analyze the Bilevel Tree-of-Hubs Location Problem with Prices (BTHLPwP). In the BTHLPwP, two hierarchical agents are involved in the decision-making process. The agent with higher hierarchy (leader) determines the design of the tree-of-hubs network (as in the THLP) and sets the prices users must pay for using the links of the hub {backbone} network. The leader's objective is to maximize profit, calculated as the total revenue from customer payments minus the maintenance and setup costs. The agent with lower hierarchy (follower) makes decisions on network usage based on the hub {backbone} network's structure and the prices set by the leader. Typically, multiple followers are involved in this decision-making process. Their objective is to minimize routing costs for their commodities.

In most HLPs, the unitary transportation costs incurred by users are typically assumed to be known and fixed by the decision maker~\citep[see e.g.][]{contreras2019hub}. In this paper, apart from making decisions about  activating hub nodes, establishing links, and routing commodities, we extend the analysis to decide also the prices paid by customers for using the hub links. This is particularly useful when constructing a (physical or virtual) hub backbone network from scratch and the investment made to construct the network want to be recovered. However, if the prices are not adequately stablished, it may have a negative effect, since if the prices are high (as optimal strategy for the decision maker maximizing its profit), the users may decide to use route the commodity by their selves or use the services of a third-party with better prices, being null the return for the leader agent. Thus, finding an equilibrium between the profit received by the leader agent and the costs incurred in the users is far from trivial, and a specialized study is required. 

We model the BTHLPwP as a Bilevel Optimization Problem. Bilevel optimization is a valuable tool for decision-making in contexts where a predefined hierarchy exists among multiple interacting agents. In these mathematical optimization problems certain constraints are defined by the solutions of another nested optimization problem. In such scenarios, one agent (the leader) makes decisions first, thereby influencing the subsequent decisions of another agent (the follower).
 
 In the BHLPwP, the leader represents the entity responsible for constructing the hub {backbone} network and setting prices, aiming to maximize profit. Conversely, the follower/s are the user/s who, based on the network structure and prices set by the leader, decide which network arcs to utilize for routing their commodities. It is assumed that, in addition to using the hub {backbone} network, each user has the option to route their commodities through an alternative delivery mode.

In this scenario, as already mentioned, the leader's profit maximization strategy might involve setting high prices for network usage. However, this could incentivize followers to opt for directly routing their commodities, resulting in reduced revenue for the leader. Conversely, setting low prices may encourage customers to use the hub {backbone} network, but the leader's profit could be insufficient to cover maintenance and setup costs. The bilevel programming approach seeks to balance the decisions of the leader and the followers. {In practice, the BTHLPwP may apply to a railway network, where the administrator sets tariffs for each segment of the network for a predetermined period, and users decide whether to use the network or opt for an alternative mode of travel or shipment. Similarly, in air cargo transportation, delivery companies determine whether to use the infrastructure network provided by another company, considering the prices for each segment used.}

\subsection{Related Works}

Hub location is an active area of research, driven by its practical applications and the mathematical challenges it presents. The core problem involves determining the optimal positions of hub nodes and links, with the objective of minimizing both setup costs and routing costs for a set of commodities. HLPs were first introduced in seminal works such as \citep{articleOKelly1,articleOKelly2}. Over the years, various problem variants have been proposed, including capacitated \citep{da2008capacitated,merakli2017capacitated,rodriguez2008solving}, uncapacitated (fixed-cost) \citep{contreras2011dynamic,topcuoglu2005solving,vasconcelos2011uncapacitated}, single and multiple allocation \citep{boland2004preprocessing,ebery2000capacitated,mohammadi2019reliable}, multi-periods \citep{khaleghi2022multi}, covering objectives \citep{peker2015p}, and those incorporating uncertainty \citep{correia2018stochastic}, among others.

In this study, we focus on a specific variant of HLPs, known as the THLP, introduced in \citep{contreras2009tight,contreras2010tree}. The THLP involves determining the locations of a specified number ($p$) of hub nodes from a given set of nodes, as well as the associated hub links, to form a tree-shaped hub {backbone} network. The objective is to minimize the total routing costs for a predefined set of commodities. This tree-shaped network configuration has proven its utility in various hub {backbone} network design contexts \citep{nickel2001hub,soylu2019multiobjective}. While the THLP was originally introduced as a single-allocation hub location problem, we adopt a multiple-allocation version, which is more practical for real-world applications.

The THLP has received increasing attention over the past decade, resulting in the development of various solution approaches. In \citep{Martins}, the authors propose a novel Benders-based branch-and-cut method for solving the problem. Another approach, utilizing a genetic algorithm, is presented in \citep{pessoa}. Additionally, several variants of the THLP have been investigated. For example, \citep{pozo} introduces a generalized aggregation method that employs different operators instead of minimizing the total routing costs. In \citep{victor}, certain hub nodes can be upgraded, leading to reduced transportation costs for routes that use those hubs. Finally, \citep{kayicsouglu2021multiple} examines a multiple-allocation version of the THLP.

While most HLPs traditionally emphasize cost-oriented objective functions, there is growing interest in exploring profit-oriented objectives. For instance, studies such as \citep{campbell1996hub} and \citep{marianov1999location} focus on maximizing some form of covered demand. Similarly, \citep{articlemax} examines a reliable version of the HLP that aims to maximize network performance.

In these studies, the goal is to design a hub backbone network that maximizes the gain derived from its usage. As a result, input parameters must include the prices paid by users for utilizing the hub {backbone} network. These prices play a critical role in influencing network usage and, consequently, its profitability. Both excessively high and low prices can negatively impact the decision maker's gain, making it essential to balance prices to optimize both user satisfaction and network manager profits. For example, \citep{luer} addresses a pricing and hub location problem in which new firms must determine strategies to manage and expand an incomplete hub {backbone} network. They propose a mixed-integer linear programming formulation and develop a genetic algorithm to solve the problem. Similarly, \cite{erdougan2022hub} investigates a hub location problem with pricing, formulating it as a non-linear bilevel programming problem. They linearize the model and enhance its performance using variable reductions and valid inequalities.

In many decision-making scenarios where a single decision maker seeks to maximize their return while being influenced by others, as in the hub {backbone} network problems discussed earlier, a hierarchical structure among decision makers is necessary to streamline the decision-making process. Such situations can arise in both cooperative and competitive systems. Bilevel optimization, as described in \citep{bard2013practical}, provides a robust framework for addressing conflicts between hierarchical decision makers within a unified mathematical programming problem. This approach enables the modeling of intricate decision-making dynamics involving multiple agents with distinct objectives and mutual influences.

Bilevel optimization has been applied to model problems related to hub {backbone} network design. For example, \citep{korani2021bi} studies a problem in which the leader designs the hub backbone network to minimize the incurred costs, while the follower determines how the commodities are shipped through the hub network, aiming to minimize the flow lost. This bilevel model considers the reliability of the designed network. A heuristic algorithm based on the KKT conditions is proposed to solve a set of realistic instances in an aviation context. In \citep{sasaki}, the authors explore a different Stackelberg HLP, where a leader firm competes with smaller firms, each aiming to maximize its profit. This problem is formulated as a bilevel optimization problem, and a sequential quadratic programming-based approach is proposed for solving it. Similarly, in \citep{cvokic1}, a leader-follower model is introduced, where each agent simultaneously constructs its hub {backbone} network to attract customers. A bilevel optimization model is formulated, followed by a MILP reformulation, and an alternating heuristic is developed for solving it. In another study by \citep{esmaeili}, a hub location problem is presented in which the leader freely selects the locations of some hub nodes, while the follower locates its hub nodes in the remaining spaces or uses hubs from the competition. This problem is modeled using bilevel optimization and solved with an imperialist competitive metaheuristic algorithm.

It is worth noting that, unlike the problems described above, our problem does not involve competition. The distinction between the two main hierarchical agents making decisions to design a hub {backbone} network arises from their differing roles in the decision-making process. The leader decides on the network's design, while the follower determines its usage. However, both decisions clearly influence each other. Furthermore, in contrast to the problems mentioned above, where the hub {backbone} network is assumed to be complete, we consider a tree-shaped hub {backbone} network in our problem. Lastly, unlike the previously mentioned works, our problem goes beyond determining the location of hub nodes and links. The leader must also decide on the prices for using the network, adding another layer of complexity to the decision-making process.

As evident from the papers mentioned above, population-based algorithms have proven successful in solving complex bilevel problems. A recent review by \citep{camacho2023metaheuristics} found that evolutionary algorithms are the most commonly used metaheuristic approaches to tackle bilevel optimization problems. However, co-evolutionary algorithms have received comparatively less attention in this domain. In co-evolutionary algorithms, the structure typically involves one population for the leader's decisions and another for the follower's decisions, with a co-evolutionary operator facilitating information exchange between the two populations \cite{camacho2023metaheuristics}. In our study, we propose an innovative approach in which we partition the leader's variables into two disjoint subsets for the populations. The followers' problem is then optimally solved in a nested manner. Specifically, we allocate one population for the pricing variables and another for the network design.

\subsection{Contributions}

In this paper, we study the BTHLPwP, and the main contributions are as follows:

\begin{enumerate}
\item We introduce, for the first time, a pricing variant of the multiple-allocation Tree-of-Hubs Location Problem, which we call the BTHLPwP.
\item We formulate the BTHLPwP as a bilevel optimization problem, where the leader is responsible for determining the activated hub nodes and links of the network, as well as setting the prices (per unit of flow) for utilizing those arcs. The followers, in turn, decide on their optimal routes for shipping commodities at minimum cost.
\item {We formulate the BTHLPwP as a bilevel mixed-integer optimization problem and provide a reformulation as a single-level mixed-integer non-linear optimization problem. We show that the non-linear terms can be linearized, yielding a mixed-integer linear reformulation of the problem.}
\item  We design a novel co-evolutionary metaheuristic algorithm specifically tailored to efficiently solve the BTHLPwP, leveraging its structural properties. This algorithm strategically integrates decisions related to tree construction, price setting, and commodity routing, offering an effective approach for solving benchmark instances within a reasonable CPU time.
\item {We conduct an extensive series of computational experiments using adapted instances from the literature to validate our proposed approach. First, the single-level reformulation is solved using Gurobi. Next, we test the proposed co-evolutionary algorithm. Finally, we perform additional experimentation by using the prices and hub backbone network obtained in the best solution from the metaheuristic as a warm start for Gurobi, aiming to optimally solve the instances and validate the effectiveness of the proposed co-evolutionary algorithm.}
\end{enumerate}

\subsection{Paper Structure}

The remainder of this paper is organized as follows: Section \ref{sec:form} presents the problem statement and the proposed bilevel programming formulation. {The equivalent single-level reformulation for the BTHLPwP is developed in Section \ref{sec:SLR}.} The co-evolutionary algorithm developed is detailed in Section \ref{sec:coea}. Section \ref{sec:exp} is devoted to computational experimentation, where the obtained results are summarized, analyzed, and discussed. Finally, the paper concludes with the key findings of this research and potential future research directions, presented in Section \ref{sec:concl}.

\section{The Bilevel Tree-of-Hubs Location Problem with Prices}\label{sec:form}

This section introduces the problem under analysis and establishes the notation used throughout the paper. Additionally, we present the proposed bilevel programming formulation for the problem.

Let $G = (I, E)$ be an {undirected} network, where $I = \{1, \ldots, m\}$ is the set of potential hub nodes and $E$ is the set of links between them (the potential inter-hub edges). For each edge $e = \{i, j\} \in E$, we denote by $f_e$ the cost of activating the edge in the backbone network provided that the hub nodes $i$ and $j$ are also activated. Without loss of generality, we assume that $G$ is a complete network, as any missing connections can be represented by fixing large activation costs for those arcs that are not available in the network.

Service demand is represented by a set of commodities defined over pairs of users in a set $J$, indexed by the set {$\mathcal{C} = \{a = (o_a, d_a, w_a) : a \in C\}$}, where the triplet $(o_a, d_a, w_a)$ indicates that an amount of flow $w_a > 0$ must be routed from origin $o_a \in J$ to destination $d_a \in J$ through the set of activated edges in $E$. {It is assumed that activating an edge $e = \{i, j\}$ in $E$ allows for bidirectional usage of the edge $\{i, j\}$, i.e., the arcs $(i, j)$ and $(j, i)$ can both be utilized. The set $A$ denotes the arcs induced by the edge set $E$, defined as $A = \{(i, j), (j, i) : \{i, j\} \in E\}$.}

Once the backbone network is established, commodities can be routed through the network with no restrictions on utilizing one or more hub arcs. Thus, the routes for the commodities {consist of sequences of arcs in the form $(o_a, k_1), (k_1, k_2), \ldots, (k_{\ell(a)-1}, k_{\ell(a)}), (k_{\ell(a)}, d_a)$, where $k_1, \ldots, k_{\ell(a)}$ are active hubs in the network, and $\{k_1, k_2\}, \{k_2, k_3\}, \ldots, \{k_{\ell(a)-1}, k_{\ell(a)}\}$ are activated inter-hub links. The \textit{access arcs}, linking the origins with the activated hubs and those linking the activated hubs with the destinations, are considered as active.}

As previously mentioned, two different agents are involved in the decisions of the BTHLPwP simultaneously. On the one hand, the leader decides which nodes and links in $G$ must be activated as {hubs and inter-hub links, respectively}. We denote by $H=(\bar{I}, \bar{E})$ the network of activated hubs and links ($\bar{I} \subseteq I$ and $\bar{E} \subseteq E$). The resulting network is assumed to have a tree structure with $p$ nodes, i.e., $\bar{H}=(\bar{I}, \bar{E})$ forms a connected acyclic graph (and consequently, with $p-1$ edges). {Tree networks are known to be useful in distribution systems where the setup and distribution costs are high, since trees are connected networks with the smallest number of direct connections between nodes. This structure helps reduce the overall infrastructure required, as well as the number of physical routes for the commodities. Thus, both capital investment in infrastructure and ongoing maintenance costs, such as repairs and upgrades, are reduced. With the above notation,} the construction of such a network incurs a cost of $\sum_{e \in \bar{E}} f_e$.

The leader also determines the prices for utilizing the constructed network. Each commodity $a \in \mathcal{C}$ that traverses the network via the arc $(i,j)$ incurs a cost (paid by the user to the leader), which is assumed to be proportional to the flow $w_a$, denoted as $w_a \pi_{ij}$ (where $\pi_{ij}$ is the price per unit flow for traversing the arc). Additionally, the maintenance of the network incurs costs. We denote by $g_i \geq 0$ the unit cost incurred by the leader when commodity $a$ uses the hub backbone network and enters node $i \in I$. {This cost (which can be set to zero) accounts for maintenance, workload, administrative paperwork, and other costs the leader incurs when a customer enters the network. We also assume that a third-party agent offers the service of routing commodity $a$ at a cost $c_a$. Given the different routing costs, the follower decides whether to route the commodity through the third-party service or the hub network.}

The leader's objective is to maximize overall profit, which consists of the total income minus the costs incurred for all commodities using the network. In contrast, once the leader has determined the structure of the hub {backbone} network, the follower routes commodities at the minimum cost, {either through this network or via the direct link offered by the third-party service.}

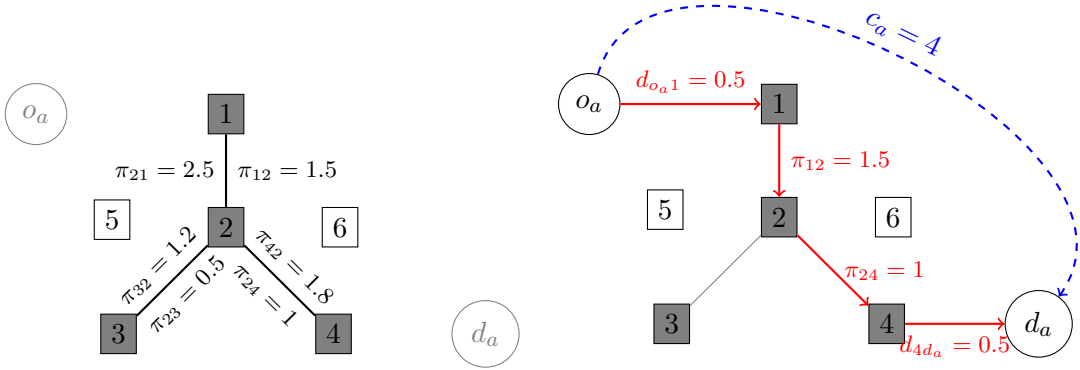
\begin{figure}[h]
\vspace*{-1cm}
{\begin{tikzpicture}[node distance=2cm, every node/.style={font=\scriptsize}, scale=0.5]

\node[draw, shape=circle, gray] (oA) at (0,10) {\color{gray} \small $o_a$};
\node[draw, fill = gray, shape=rectangle] (H1) at (5,10) {\small $1$};
\node[draw, fill = gray, shape=rectangle] (H2) at (5,7) {\small $2$};
\node[draw,  shape=rectangle] (H5) at (2,7.2) {\small $5$};
\node[draw,  shape=rectangle] (H6) at (8,7) {\small $6$};
\node[draw, fill = gray, shape=rectangle, below left of=H2] (H3) {\small $3$};
\node[draw, fill = gray, shape=rectangle, below right of=H2] (H4) {\small $4$};
\node[draw, shape=circle, right of=H4, gray] (dA) {\color{gray}\small $d_a$};

\draw[thick] (H1) -- (H2) node[midway, right] {$\pi_{12}=1.5$} node[midway, left] {$\pi_{21}=2.5$};
\draw[thick] (H2) -- (H3)  node[midway,below,sloped] {$\pi_{23}=0.5$} node[midway,above,sloped] {$\pi_{32}=1.2$};
\draw[thick] (H2) -- (H4) node[midway,below,sloped]{$\pi_{24}=1$} node[midway,above,sloped] {$\pi_{42}=1.8$};

\draw[->, bend left=100, dashed, white] (oA) to node[right] {\small \color{white}$c_a=6$} (dA);

\end{tikzpicture}}~\hspace*{-1cm}{\begin{tikzpicture}[node distance=2cm, every node/.style={font=\scriptsize}, scale=0.5]

\node[draw, shape=circle] (oA) at (0,10) { \small $o_a$};
\node[draw, fill = gray, shape=rectangle] (H1) at (5,10) {\small $1$};
\node[draw, fill = gray, shape=rectangle] (H2) at (5,7) {\small $2$};
\node[draw,  shape=rectangle] (H5) at (2,7.2) {\small $5$};
\node[draw,  shape=rectangle] (H6) at (8,7) {\small $6$};
\node[draw, fill = gray, shape=rectangle, below left of=H2] (H3) {\small $3$};
\node[draw, fill = gray, shape=rectangle, below right of=H2] (H4) {\small $4$};
\node[draw, shape=circle, right of=H4] (dA) {\small $d_a$};

\draw[->, thick, red] (oA) -- (H1) node[midway, above] {$d_{o_{a}1}=0.5$};
\draw[->,  thick, red] (H1) -- (H2) node[midway, right] {$\pi_{12}=1.5$};
\draw[gray] (H2) -- (H3);
\draw[->,  thick, red] (H2) -- (H4) node[midway, right]{$\pi_{24}=1$};
\draw[->, thick, red] (H4) -- (dA) node[midway, below] {$d_{4d_{a}}=0.5$};

\draw[->, bend left=100, dashed, blue, thick] (oA) to node[above, sloped] {\color{blue}\small $c_a=4$} (dA);

\end{tikzpicture}}
\caption{Illustrative example showcasing the decision-making process in the BTHLPwP. \label{fig:0}}
\end{figure}

{In Figure \ref{fig:0}, we illustrate the different phases of the problem under analysis. In the left picture, we show a feasible solution for the leader when $p=4$ (i.e., 4 hub nodes to activate). A set of 6 potential locations for the hub nodes is given (labeled as $1$, $2$, $\ldots$, $6$), represented in the plot with squares. Among them, the leader decides to activate nodes $1$, $2$, $3$, and $4$ (highlighted in gray). The tree structure is also decided, with the inter-hub links being $\{1,2\}$, $\{2,3\}$, and $\{2,4\}$. The prices for these links are set to $\pi_{12}=1.5$, $\pi_{23}=0.5$, and $\pi_{24}=1$.}

{With the leader's decisions already set, the followers (customers) decide how to route the commodities. In the right plot, we show a possible route for a commodity with origin $o_a$, destination $d_a$, and flow $w_a=1$. On the one hand, the follower decides whether to use the constructed hub network or the direct/third-party distribution system. If the commodity is routed through the direct link (dashed blue arrow), the overall transportation cost is $w_a c_a = 1 \cdot 4 = 4$. On the other hand, if the commodity is routed via the hub network, the shortest path linking $o_a$ and $d_a$ using the access links and the inter-hub links is chosen. In this case, the shortest route is highlighted in red, being the path formed by the arcs $(o_a, 1)$, $(1,2)$, $(2,4)$, and $(4,d_a)$, with the overall cost being $1 \cdot (0.5 + 1.5 + 1 + 0.5) = 3.5$. Thus, since routing the commodity through the hub network is less costly than using the direct route offered by the third-party service, the commodity would be routed through the hub network. The followers incur the routing costs (3.5 units for this particular commodity), whereas the leader obtains as profit the price for using the hub network (those incurred by the follower) minus the cost of activating the inter-hub links and the maintenance costs for those commodities that use the hub network.}

{The leader's and followers' decisions are detailed below using the notation considered in our proposed mathematical model.}

\subsection{Leader's Decisions} 

The leader's decisions are captured through the following sets of variables, corresponding to their three main responsibilities: hub activation, inter-hub link activation, and pricing.
$$
y_i = \left\{\begin{array}{cl} 1 & \mbox{if hub $i$ is activated.}\\
0 & \mbox{otherwise}
\end{array}\right. \quad \forall i \in I,
$$
$$
z_{ij}  = \left\{\begin{array}{cl} 1 & \mbox{if edge $\{i,j\}$ is activated}\\
0 & \mbox{otherwise}
\end{array}\right. \quad \forall \{i, j\}\in E,
$$
$$
\pi_{ij} \geq 0: \mbox{price for routing a unit flow through arc $(i,j)$}, \forall i, j \in I.
$$

\subsection{Followers' Decisions}

Given the hub {backbone network ($\bar{H}$), and the prices ($\pi$) set by the leader, the follower(s) determine(s) the utilization of the hub backbone network}. For a commodity $a \in \mathcal{C}$, the decision involves whether to route it via the hub {backbone} network using its hub arcs.  {It is important to note that a user (follower) may have multiple commodities to route. However, since no capacity constraints are considered, the routing costs for all commodities can be summed into a single term that accounts for the costs of all commodities. Therefore, the case of multiple commodities with independent followers can be reduced to a single follower scenario. This is supported by some of the results presented in \citep{calvete2007linear}. }

The cost of utilizing arc $(i,j)$, denoted by $C_{a,i(a),j(a)}$, includes the collection, distribution, and transportation costs for commodity $a$ and can be expressed as:
$$
C_{a,i(a),j(a)} = w_a \Big( \mathrm{d}_{o_a,i} + \sum_{s=1}^{\ell(a)-1} \pi_{k_s k_{s+1}} + \mathrm{d}_{k_{\ell(a)},o_a}\Big)
$$
{where $\mathrm{d}_{o_a,i}$ represents the unit collection cost to the first node when commodity $a$ enters the hub {backbone} network. The routing sequence within the hub {backbone} network is given by the arcs $(k_1, k_2), \ldots, (k_{\ell(a)-1}, k_{\ell(a)})$. Finally, $\mathrm{d}_{\ell(a),d_a}$ denotes the unit distribution cost from the last hub node in the path to the destination of commodity $a$, denoted as $d_a$.}

If the follower opts to route commodity $a$ directly through a third-party service, bypassing the hub {backbone} network, a cost of $c_a$ is incurred. The follower is assumed to select the least costly option for routing each commodity. {Consequently, the follower's objective is to minimize the total routing costs by determining whether to use the hub {backbone} network established by the leader or to route the commodity directly without utilizing the hub {backbone} network.} In the first case, utilizing the hub network generates profit for the leader, whereas in the second scenario, the leader receives no income. Note that, since different commodities are involved in the follower's decisions, the routing for each commodity is determined independently.

The routes for the commodities are modeled in our optimization problem using the following sets of variables:

$
x_{aij} =\left\{\begin{array}{cl} 1 & \mbox{If commodity $a$ uses the {inter-hub}}\\
 & \mbox{arc $(i,j)$ in its route.}\\
0 & \mbox{otherwise}
\end{array}\right., \forall a \in \mathcal{C}, i, j \in I.
$

$
x^1_{ai}=\left\{\begin{array}{cl} 1 & \mbox{If commodity $a$ is initially routed}\\
& \mbox{via hub node $i$.}\\
0 & \mbox{otherwise}
\end{array}\right. \quad \forall a \in \mathcal{C}, i \in I.
$

$
x^2_{aj}=\left\{\begin{array}{cl} 1 & \mbox{If commodity $a$ is endly routed}\\
& \mbox{via hub node $j$.}\\
0 & \mbox{otherwise}
\end{array}\right. \quad \forall a \in \mathcal{C}, j \in I.
$

$
q_a=\left\{\begin{array}{cl} 1 & \mbox{If commodity $a$ is routed directly}\\
& \mbox{from origin $o_a$ to destination $d_a$.}\\
0 & \mbox{otherwise}
\end{array}\right. \quad \forall a \in \mathcal{C}.
$

\subsection{A Bilevel Optimization Model for the BTHLPwP}

All the sets, parameters, and decision variables introduced above are summarized in Table \ref{table:notation}.

\begin{table}[H]
\begin{tabular}{|cp{11cm}|}
\multicolumn{2}{c}{\bf Index Sets}\\\hline
$I$ & Set of potential hub nodes,  where $|I|=m$.\\
$\mathcal{C}$ & Set of commodities $a \equiv (o_a, d_a, w_a)$.\\\hline
\multicolumn{2}{c}{\bf Parameters}\\\hline
$o_a$ & Origin node of commodity $a\in \mathcal{C}$.\\
$d_a$ & Destination node of commodity $a\in \mathcal{C}$.\\
$w_a$ & Flow of commodity $a\in \mathcal{C}$ to be routed.\\
$g_i$ & Maintenance cost for entering the hub backbone network through node $i\in I$.\\
$f_{ij}$ & Set-up cost for inter-hub link $\{i,j\}$, $i, j\in I$.\\
$\mathrm{d}_{o_a,i}$ & Unit collection cost for routing commodity $a\in \mathcal{C}$ through the hub backbone network when entering at the hub node $i\in I$.\\
$\mathrm{d}_{i,d_a}$ & Unit collection cost for routing commodity $a\in \mathcal{C}$ through the hub backbone network when leaving from the hub node $i\in I$.\\
$c_a$ & Cost for routing commodity $a\in \mathcal{C}$ by the third-party agent without using the hub backbone network.\\\hline
\multicolumn{2}{c}{\bf Decision Variables}\\\hline
$y_i \in \{0,1\}$ & Activation of hub node $i\in I$.\\
$z_{ij} \in \{0,1\}$ & Activation of hub link $\{i,j\}$, $i,j \in I$.\\
$\pi_{ij} \geq 0$ & Price for routing a unit flow through arc $(i,j)$.\\\hline
$x_{aij} \in \{0,1\}$ & Use of arc $(i,j)$ when routing commodity $a\in \mathcal{C}$.\\
$x_{ai}^1 \in \{0,1\}$ & Entering the hub backbone network through node $i\in I$ for commodity $a\in \mathcal{C}$.\\
$x_{aj}^2 \in \{0,1\}$ & Leaving the hub backbone network through node $j\in I$ for commodity $a\in \mathcal{C}$.\\
$q_a \in \{0, 1\}$ & Use the third-party service for routing commodity $a \in \mathcal{C}$ without using the hub {backbone} network.\\\hline
\end{tabular}
\caption{Index sets, parameters, and variables used in our mathematical optimization model.\label{table:notation}}
\end{table}

Using the notation mentioned above, we propose the following bilevel program to model the BTHLPwP:

\begin{align}
\max_{y, z, \pi} && &\dsum_{i, j \in I}  \pi_{ij} \dsum_{a \in \mathcal{C}} w_{a} x_{aij} - \dsum_{a \in \mathcal{C}}\dsum_{i  \in I} w_{a}g_{i} x^{1}_{ai} - \dsum_{i, j  \in I} f_{ij} z_{ij} \label{l:0}\\
\mbox{s.t. } & &&\dsum_{i \in I} y_i =p,\label{l:1}\\
&& &\dsum_{i, j \in I:\atop i>j} z_{ij} = p-1,\label{l:2}\\
&& &z_{ij} \leq y_i, \qquad \forall i, j \in I \label{l:3}\\
&& &z_{ij} \leq y_{j}, \qquad \forall i, j \in I \label{l:4}\\
&& &\dsum_{j \in I \backslash S } z_{ij} \geq 1, \qquad \forall S \subseteq I, i \in S \label{l:5}\\
&& &y_i \in \{0,1\}, \quad \quad \ \ \forall i \in I \label{l:6}\\
&& &z_{ij} \in \{0,1\}, \qquad \ \forall i, j \in I \label{l:7}\\
&& &\pi_{ij} \geq 0, \qquad \qquad \forall i, j \in I \label{l:8}\\
&& & \text{where for a leader's solution } (y,z,\pi), \ (x,x^1,x^2,q) \text{ solve: } \nonumber \\
& &\min_{x,x^1,x^2,q}& {\dsum_{a \in \mathcal{C}} w_{a} \Big( \dsum_{i\in I}  (\mathrm{d}_{o_a,i} x^1_{ai}  + \mathrm{d}_{i,d_a}x^2_{ai})  +\dsum_{i, j \in I} \pi_{ij} x_{aij} + c_{a} q_{a}\Big)}\label{f:0}\\
&&\mbox{s.t. } & \dsum_{i \in I} x^1_{ai} + q_{a} =1, \qquad \forall a \in \mathcal{C} \label{f:1}\\
&	&	& \dsum_{i \in I} x^2_{ai} + q_{a} =1,\qquad \forall a \in \mathcal{C} \label{f:2}\\
&	&	& \dsum_{j \in I} x_{aij} + x^1_{ai}  - \dsum_{j \in I} x_{aj i} - x^2_{ai}  =0, \qquad \forall a \in \mathcal{C}, i \in I \label{f:3}\\
&	&	& x_{ai}^1 \leq y_i,\qquad \forall a \in \mathcal{C}, i \in I  \label{f:4}\\
&	&	&x_{ai}^2 \leq y_i,\qquad \forall a \in \mathcal{C}, i \in I \label{f:5}\\
&	&	& x_{aij}+x_{aji} \leq z_{ij},\qquad \forall a \in \mathcal{C}, i, j \in I, (i<j)\label{f:6}\\
&	&	&  x_{aij} \in \{0,1\}, \qquad \forall a \in \mathcal{C}, i,j \in I\label{f:7}\\
&	&	& x^1_{ai} \in \{0,1\},\qquad \ \forall a \in \mathcal{C}, i \in I\label{f:8}\\
&	&	&  x^2_{ai} \in \{0,1\},\qquad \ \forall a \in \mathcal{C}, i \in I\label{f:9}\\
&	&	&  q_{a} \in \{0,1\},\qquad \ \ \forall a \in \mathcal{C}.\label{f:10}
\end{align}

The leader's objective function \eqref{l:0} maximizes the overall profit from network usage, which is the difference between incomes and costs. Constraints \eqref{l:1} and \eqref{l:2} ensure that exactly $p$ hubs and $p-1$ hub links are activated. Constraints \eqref{l:3} and \eqref{l:4} prevent the activation of hub links unless both endpoints are activated. Constraints \eqref{l:5} enforce the connectivity of the hub {backbone} network, which, in conjunction with the previous constraints, guarantees that the resulting hub link network forms a tree. Finally, constraints \eqref{l:6}, \eqref{l:7}, and \eqref{l:8} define the domain for the leader's  decision variables.

On the other hand, constraints \eqref{f:0} to \eqref{f:10} define the follower's problem, which is parameterized by the leader's variables (activated hubs $y$, activated links $z$, and prices $\pi$).  The objective function \eqref{f:0} minimizes the total routing cost for all commodities. Constraints \eqref{f:1} and \eqref{f:2} ensure that either a unique hub node is used for the initial and final routing of the commodity, or the commodity is routed directly without using the hub {backbone} network. Constraints \eqref{f:3} represent the flow conservation constraints, ensuring that each commodity is properly routed if it chooses to use the hub {backbone} network. Constraints \eqref{f:4}, \eqref{f:5}, and \eqref{f:6} prevent routing commodities through non-activated hubs or links. Finally, constraints \eqref{f:7}, \eqref{f:8}, \eqref{f:9}, and \eqref{f:10} enforce the binary nature of the follower's decision variables.

Note that if the hub {backbone} network (variables $y$ and $z$) and the prices (variables $\pi$) are known, the follower's problem simplifies to a minimum cost flow problem with unitary flows. This problem can be efficiently solved (in polynomial time) by solving $|\mathcal{C}|$ separate shortest path problems using Dijkstra's algorithm for each commodity. If the cost of the shortest path for a commodity $a$ is smaller than $c_a$, the commodity is routed through the hub {backbone} network; otherwise, the commodity is routed directly via the third-party service, i.e., $q_a=1$. An assumption must be made when multiple shortest paths are obtained for a specific commodity given the leader's decision. In this case, the optimistic approach is assumed \citep{kalashnikov2015bilevel}, \citep{dempe2002foundations}. Specifically, the shortest path that yields the highest profit for the leader is chosen as the follower's decision.

{

\section{An Equivalent Single-Level Reformulation for the BTHLPwP}\label{sec:SLR}

As mentioned in the previous section, when the leader's decision variables $y$ and $z$ associated with the hub backbone network, and the pricing variables $\pi$, are known, the parameterized follower's problem reduces to a minimum cost flow problem. Consequently, the bilevel problem defined by Eqs. (\ref{l:0})–(\ref{f:10}) can be reformulated as an equivalent single-level program using duality theory and the optimality conditions of the follower's problem. To achieve this reformulation, we consider the case where the leader has fixed the prices ($\bar{\pi}$) and constructed the hub backbone network ($\bar{H}, \bar{T}$). The resulting follower's problem is as follows:

\begin{align}
\min_{x,x^1,x^2,q}& \dsum_{a \in \mathcal{C}} w_{a} \Big( \dsum_{i \in \bar{H}:\atop \{i,j\} \in \bar{T}} (\mathrm{d}_{o_a,i} x^1_{ai}  + \mathrm{d}_{i,d_a}x^2_{ai})  +\dsum_{i,j \in \bar{H}:\atop \{i,j\} \in \bar{T}} \bar{\pi}_{ij} x_{aij} + c_{a} q_{a}\Big)\label{fparam:0}\\
\mbox{s.t. } & \dsum_{i \in \bar{H}} x^1_{ai} + q_{a} =1,\qquad  \forall a \in \mathcal{C}, \label{fparam:1}\\
		& \dsum_{i \in \bar{H}} x^2_{ai} + q_{a} =1, \qquad \forall a \in \mathcal{C},\label{fparam:2}\\
		& \dsum_{j \in \bar{H}:\atop \{i,j\} \in \bar{T}} x_{aij} + x^1_{ai}  - \dsum_{j \in \bar{H}:\atop \{i,j\} \in \bar{T}} x_{aj i} - x^2_{ai}  =0, \qquad \forall a \in \mathcal{C}, i \in \bar{H}\label{fparam:3}\\
		&  x_{aij} \in \{0,1\}, \qquad \forall a \in \mathcal{C}, i,j \in \bar{H}: \{i,j\} \in \bar{T} \label{fparam:7}\\
		& x^1_{ai} \in \{0,1\},\qquad \ \forall a \in \mathcal{C}, i \in \bar{H} \label{fparam:8}\\
		&  x^2_{ai} \in \{0,1\},\qquad \ \forall a \in \mathcal{C}, i \in \bar{H} \label{fparam:9}\\
		&  q_{a} \in \{0,1\},\qquad \ \ \forall a \in \mathcal{C}.\label{fparam:10}
	\end{align}

Two key aspects of the follower's problem defined in Eqs. \eqref{fparam:7}-\eqref{fparam:10}, are worth highlighting. First, since the hub backbone network is known, constraints \eqref{f:4}-\eqref{f:6} are no longer necessary. Second, the decision variables can be relaxed to continuous variables, i.e., on can replace constraints \eqref{fparam:7}-\eqref{fparam:10} by:
\begin{align*}
    		&  x_{aij} \geq 0, \qquad \forall a \in \mathcal{C}, i,j \in \bar{H}: \{i,j\} \in \bar{T},\\
		& x^1_{ai} \geq 0,\qquad \ \forall a \in \mathcal{C}, i \in \bar{H},\\
		&  x^2_{ai} \geq 0,\qquad \ \forall a \in \mathcal{C}, i \in \bar{H},\\
		&  q_{a} \geq 0,\qquad \ \ \forall a \in \mathcal{C}.
\end{align*}
Then, the resulting follower problem (when the leader decision variables are fixed) is a linear optimization problem, that can be separated for each commodity $a\in \mathcal{C}$. Therefore, its dual can be constructed in a straightforward manner for each $a \in \mathcal{C}$. Let $\lambda^1_a$, $\lambda^2_a$, $\mu_a$ be the dual variables associated to the functional constraints \eqref{fparam:1}, \eqref{fparam:2}, and  \eqref{fparam:3}, respectively. 
Therefore, the dual model for the parameterized follower's problem for commodity $a \in \mathcal{C}$ reads:

\begin{align}
\max & \;\; \lambda_{a}^1 + \lambda_a^2  \label{dual:0}\\
\mbox{s.t. }          & \lambda_a^1+\lambda_a^2 \leq w_ac_a,  \label{dual:4}\\
& \lambda_a^1 - \mu_{ai} - v^1_{aij} \leq w_a \mathrm{d}_{o_a,i},\qquad \forall  i \in  \bar{H} \label{dual:1}\\
		& \lambda_a^2 + \mu_{ai} \leq w_a \mathrm{d}_{i,d_a},\qquad \forall i \in  \bar{H} \label{dual:2}\\
		& \mu_{ai}- \mu_{aj}  \leq w_a \bar{\pi}_{ij}, \qquad \forall i,j \in \bar{H}: \{i,j\} \in \bar{T} \label{dual:3a}\\
        & \mu_{aj}- \mu_{ai}   \leq w_a \bar{\pi}_{ij}, \qquad \forall i,j \in \bar{H}: \{i,j\} \in \bar{T} \label{dual:3b}\\
		&  \lambda_a^1, \lambda_a^2, \mu_{ai} \in \mathbb{R}, \qquad \forall i \in \bar{H} \label{dual:5}
\end{align}

On the other hand, if constraints \eqref{f:4}--\eqref{f:6} are introduced as capacity constraints (upper bounded by the value of the $z$ and $y$ leader variables) in the minimum cost flow problem instead of using $\bar H$ and $\bar T$, then dual variables $u_{ai}^1$, $u_{ai}^2$, and $v_{aij}$ are associated to these constraints, respectively. Then, replacing the notation for the open hubs, $\bar H$, and open links, $\bar T$, by the variables determining their structure, $y$ and $z$, respectively, we get the following formulation. the sets  Hence, the following dual formulation is obtained:
\begin{align}
\max & \;\; \lambda_{a}^1 + \lambda_a^2  - \dsum_{i \in I} (v^1_{ai}+v^2_{ai}) y_i - \dsum_{i, j \in I: i<j} u_{aij} z_{ij}\label{dual:0}\\
\mbox{s.t. }          & \lambda_a^1+\lambda_a^2  \leq w_ac_a, \label{dual:4}\\
& \lambda_a^1 - \mu_{ai} - v^1_{aij} \leq w_a \mathrm{d}_{o_a,i},\qquad \forall  i \in  I, \label{dual:1}\\
		& \lambda_a^2 + \mu_{ai} - v^2_{aij} \leq w_a \mathrm{d}_{i,d_a},\qquad \forall i \in  I,  \label{dual:2}\\
		& \mu_{ai}- \mu_{aj}  - u_{aij} \leq w_a {\pi}_{ij}, \qquad \forall  i<j \in I, \label{dual:3a}\\
        & \mu_{aj}- \mu_{ai}  - u_{aij} \leq w_a {\pi}_{ij}, \qquad \forall i<j \in I, \label{dual:3b}\\
		&  \lambda_a^1, \lambda_a^2, \mu_{ai} \in \mathbb{R}, \qquad \forall  i \in I, \label{dual:5}\\
        & v^1_{ai}, v^2_{ai}  \geq 0, \qquad \forall i \in I,\label{dual:5b}\\
        & u_{aij} \geq 0, \qquad \forall i<j \in I. \label{dual:6}
\end{align}

Hence, we can insert the above problem into the leader decision problem, and then the bilevel problem proposed for the BTHLPwP can be reformulated as the following single-level mixed-integer non-linear problem:
\begin{align}
\max && &\dsum_{i, j \in I}  \pi_{ij} \dsum_{a \in \mathcal{C}} w_{a} x_{aij} - \dsum_{a \in \mathcal{C}}\dsum_{i  \in I} w_{a}g_{i} x^{1}_{ai} - \dsum_{i, j  \in I} f_{ij} z_{ij} \label{slrNL:0}\\
\mbox{s.t. } & &&\dsum_{i \in I} y_i =p,\label{slrNL:1}\\
&& &\dsum_{i, j \in I:\atop i>j} z_{ij} = p-1,\label{slrNL:2}\\
&& &z_{ij} \leq y_i, \qquad \forall i, j \in I \label{slrNL:3}\\
&& &z_{ij} \leq y_{j}, \qquad \forall i, j \in I \label{slrNL:4}\\
&& &\dsum_{j \in I \backslash S } z_{ij} \geq 1, \qquad \forall S \subseteq I, i \in S \label{slrNL:5}\\
&& & \dsum_{i \in I} x^1_{ai} + q_{a} =1, \qquad \forall a \in \mathcal{C} \label{slrNL:6}\\
&	&	& \dsum_{i \in I} x^2_{ai} + q_{a} =1,\qquad \forall a \in \mathcal{C} \label{slrNL:7}\\
&	&	& \dsum_{j \in I} x_{aij} + x^1_{ai}  - \dsum_{j \in I} x_{aj i} - x^2_{ai}  =0, \qquad \forall a \in \mathcal{C}, i \in I \label{slrNL:8}\\
&	&	& x_{ai}^1 \leq y_i,\qquad \forall a \in \mathcal{C}, i \in I  \label{slrNL:9}\\
&	&	&x_{ai}^2 \leq y_i,\qquad \forall a \in \mathcal{C}, i \in I \label{slrNL:10}\\
&	&	&x_{aij}+x_{aji} \leq z_{ij},\qquad \forall a \in \mathcal{C}, i, j \in I, (i<j)\label{slrNL:11}\\
&&& \lambda_a^1+\lambda_a^2  \leq w_ac_a, \qquad \forall a \in \mathcal{C}\label{slrNL:12}\\
&&& \lambda_a^1 - \mu_{ai} - v^1_{aij} \leq w_a \mathrm{d}_{o_a,i},\qquad \forall a \in \mathcal{C}, i \in  I, \label{slrNL:13}\\
&&& \lambda_a^2 + \mu_{ai} - v^2_{aij} \leq w_a \mathrm{d}_{i,d_a},\qquad \forall a \in \mathcal{C}, i \in  I,  \label{slrNL:14}\\
&&& \mu_{ai}- \mu_{aj}  - u_{aij} \leq w_a {\pi}_{ij}, \qquad \forall a \in \mathcal{C}, i<j \in I, \label{slrNL:15}\\
&&& \mu_{aj}- \mu_{ai}  - u_{aij} \leq w_a {\pi}_{ij}, \qquad \forall a \in \mathcal{C}, i<j \in I, \label{slrNL:16}\\
& & & w_{a} \Big( \dsum_{i\in I} (\mathrm{d}_{o_a,i} x^1_{ai}  + \mathrm{d}_{i,d_a}x^2_{ai})  +\dsum_{i, j \in I} \pi_{ij} x_{aij} + c_{a} q_{a}\Big) = \nonumber \\
& & &  (\lambda_{a}^1 + \lambda_a^2) + \dsum_{i \in I} (v_{ai}^1+v_{ai}^2)y_i +\dsum_{i,j \in I \atop i<j} u_{aij}z_{ij}, \forall a \in \mathcal{C}\label{slrNL:16}\\
&& &y_i, z_{ij} \in \{0,1\}, \quad \quad \ \ \forall i,j \in I \label{slrNL:17}\\
&& &\pi_{ij} \geq 0, \qquad \qquad \forall i, j \in I \label{slrNL:18}\\
&	&	&  x_{aij}, x^1_{ai}, x^2_{ai}, q_{a} \in \{0,1\}, \qquad \forall a \in \mathcal{C}, i,j \in I\label{slrNL:19}\\
&&&  \lambda_a^1, \lambda_a^2, \mu_{ai} \in \mathbb{R}, \qquad \forall  i \in I, \label{slrNL:17}\\
&&& v^1_{ai}, v^2_{ai}  \geq 0, \qquad \forall i \in I,\label{slrNL:18}\\
&&& u_{aij} \geq 0, \qquad \forall i<j \in I. \label{slrNL:19}
\end{align}

In the reformulation presented in Eqs. \eqref{slrNL:0}-\eqref{slrNL:19}, Eq. \eqref{slrNL:16} ensures that the objective functions of both the primal and dual follower problems are equal, which occurs only at the optimal solutions. Consequently, it is guaranteed that the follower's problem is solved optimally. It is important to emphasize that this reformulation includes non-linear terms (see Eqs. \eqref{slrNL:0} and \eqref{slrNL:16}). Specifically, the terms $y_i (v_{ai}^1 + v_{ai}^2)$, $z_{ij}v_{aij}$, and $\pi_{ij}x_{aij}$ involve variable multiplication of binary variables by continuous ones. These expressions, however, can be linearized by introducing the following auxiliary variables:
\begin{center}
$\xi_{ai} = y_i (v_{ai}^1 + v_{ai}), \quad  \rho_{aij} = z_{ij}v_{aij}, \quad \text{and} \quad \theta_{aij} = \pi_{ij}x_{aij}$ 
\end{center}
and applying the  appropriate McCormick envelopes

Although an equivalent single-level reformulation can be developed for the optimization model above using the standard KKT optimality conditions, solving the problem exactly poses a significant computational challenge. The presence of non-linear inequalities, along with the need to consider cutset inequalities \eqref{l:5}, precludes the use of commercial solvers for solving real-world instances. While it is possible to design a relax-and-fix approach (via callbacks in off-the-shelf optimization software) to avoid incorporating the exponentially many constraints and include them as needed within a branch-and-cut approach, the number of cuts and iterations is computationally prohibitive even for small-size instances. Furthermore, the single-level reformulation requires incorporating big-M constants to linearize some of the non-linear (and non-convex) expressions. As observed in different works \citep{kleinert2020there, kleinert2023there, pineda2019solving}, determining the appropriate big-M constants is as challenging as solving the bilevel problem itself.

Given the practical limitations of the aforementioned reformulation, the following section focuses on proposing a novel metaheuristic algorithm to obtain high-quality solutions within reasonable computational time. As demonstrated in our computational experiments, these solutions ara particularly usefuls as warm-start solutions for solving, exactly, the single-level reformulation.
}

\section{{A Novel Co-Evolutionary Algorithm for the BTHLPwP}}\label{sec:coea}

In this section, the proposed Co-Evolutionary Algorithm (Co-EA) for solving the BTHLPwP is presented. Co-EAs are a specific class of Evolutionary Algorithms (EAs) used for heuristically solving optimization problems, inspired by the evolution of species.

In EAs, a population of individuals is defined based on the problem, and they undergo evolutionary operators aimed at improving their fitness. These operators include crossover, where individuals are combined to generate offspring, and mutation, where offspring may undergo small changes to add diversity to the population. A survival criterion is applied, and the process is iterated over generations. The ultimate goal of EA operators is to find solutions that enhance the fitness of the population while maintaining diversity among individuals. Therefore, the main components of EAs, which are detailed in this section for the BTHLPwP, include the fitness measure for individuals (quality of solutions), the selection procedure (maintained solutions), crossover (combination of solutions), and mutation (diversification of solutions). EAs have demonstrated significant advantages in efficiently solving computationally costly problems \citep{sloss20202019}, \citep{davis2012evolutionary}, \citep{yu2010introduction}.

In Co-EAs, two or more populations are considered simultaneously, and evolution occurs through interactions among individuals from different populations. Two main types of Co-EAs exist based on the interaction among individuals: cooperative and competitive. In cooperative Co-EAs, both populations evolve by exchanging information to facilitate their improvement, whereas in competitive Co-EAs, individuals of one population evolve at the expense of the other. For further details on Co-EAs, interested readers are referred to \citep{potter1994cooperative}, \citep{popovici2012coevolutionary}, \citep{jan}, and the references therein.

Co-EAs have proven successful in solving bilevel programming problems, attributed to the partitioning of leader and follower decision variables \citep{articleCO1}, \citep{articleCO2}, \citep{articleCO3}. In such problems, a population is created for each decision-making level (leader and follower), and both populations evolve while exchanging information between them. However, the algorithm we propose does not adhere to this conventional approach. {Instead, the two populations considered are associated with the leader's decisions, and the follower's problem is optimally solved when necessary; this is the so-called nested approach.} Thus, co-evolution between populations occurs solely in the leader's decision space.

As previously mentioned, the leader's decisions encompass two aspects: (1) constructing the hub backbone network under a tree-shaped topology, and (2) determining the unit flow prices for utilizing the network. To each decision, we assign its own population. Therefore, the first population consists of trees formed by the selected hubs, while the second comprises vectors representing the prices imposed for network usage. An asynchronous exchange of information is facilitated among these populations. In detail, a random population of trees is generated to initialize the algorithm. Subsequently, based on the generated trees, the follower's problem is optimally solved, thereby determining acceptable lower and upper bounds for the optimal prices. Then, random pricing vectors are generated within the computed intervals, with the extremes representing the lower and upper bounds. Following this, the fitness of an individual is evaluated, allowing both populations to evolve independently. When the prices are updated, co-evolution is triggered to revise the acceptable intervals for pricing. This iterative process continues until a predetermined stopping criterion is met.

\begin{figure}[H]
\centering
\includegraphics[scale=0.7]{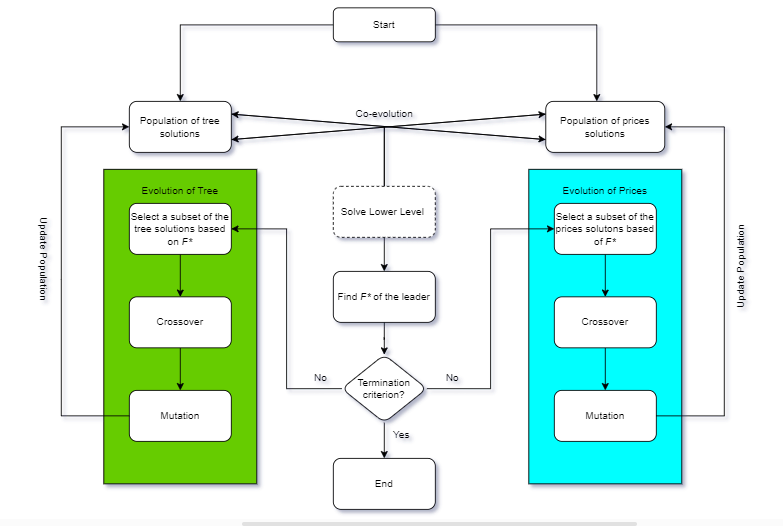}
\caption{Flowchart for the proposed Co-EA}
\label{diagramadeflujogral}
\end{figure}

A flowchart of the proposed Co-EA is shown in Fig. \ref{diagramadeflujogral}.

Next, we will detail each of the components involved in the proposed Co-EA.

\subsection{Constructing Bilevel Feasible Solutions} \label{sec:constructing}

As already mentioned, the Co-EA starts by creating two populations. One population is related to the hub backbone network, which is formed by the activated hubs and connection links under a tree topology. The other population is based on the prices imposed by the leader. Matching one solution from each population will lead to a complete leader's solution, and the follower's problem must be solved to evaluate the fitness. It is worth highlighting that both populations evolve independently, but there is an interchange of information controlled by the co-evolutionary operator.

\subsubsection{Creation of Tree-Structured Networks}
First, tree-structured hub backbone networks for selected subsets of hubs are created. These decisions are represented by the decision variables $y$ and $z$ in our mathematical model. Recall that given a set of potential sites for activating a hub $I$, exactly $p$ of them must be activated. This is represented by a binary vector $y$ of length $|I|$, where $1$ indicates that the hub is activated and $0$ otherwise.

\begin{figure}[H]
    \centering
    \includegraphics[scale = 0.5]{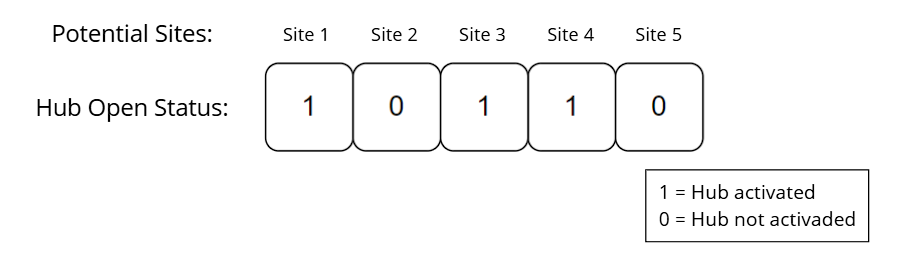}
    \caption{Example with 5 potential sites and 3 hubs to be activated.}
    \label{fig:CodeH}
\end{figure}

Based on the activated hubs, random trees on these nodes ($z$) are built, encoded using their adjacency matrices. It is worth mentioning that these matrices are symmetric, as the flow routed through the hubs is bidirectional. A representation of this partial leader's solution is depicted in Fig. \ref{fig:CodeH}.
\subsubsection{Creation of Pricing Vectors}
\begin{figure}[H]
    \centering
    \includegraphics[scale = 0.35]{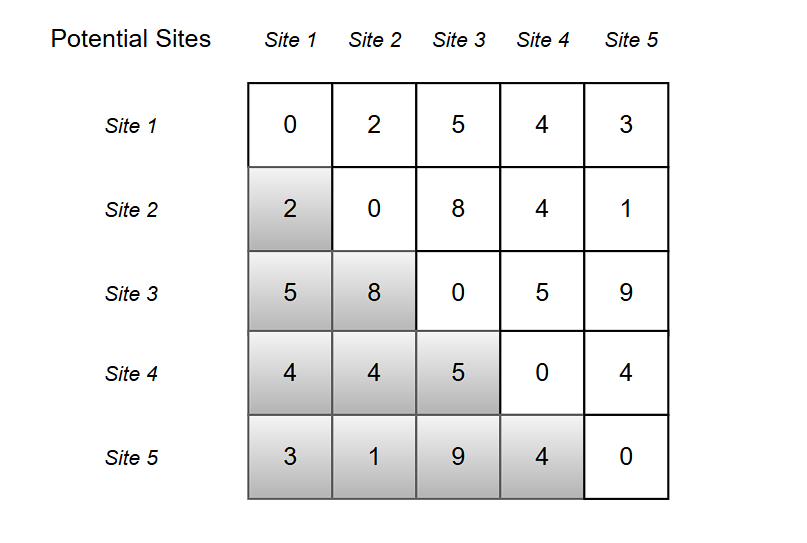}
    \caption{Example of a pricing solution: Each entry represents the price of using a specific arc. For instance, if the follower uses the arc (3,2), then a cost of 8 is incurred.}
    \label{fig:CodeP}
\end{figure}

For creating a random population of pricing vectors ($\pi$), the trees induced by $y$ and $z$ are used. Since $\pi \geq 0$, the lower bound is trivial. On the other hand, an upper bound for each edge in the tree is calculated, and acceptable intervals to generate the prices are then identified. The reasoning behind this idea is to avoid wasting effort by setting excessive prices that will discourage customers from using the network's infrastructure. That is, the commodities will be routed directly from origin to destination without returning any gain to the leader. A depiction of these prices is shown through a two-dimensional array, as in Figure \ref{fig:CodeP}.

Since the prices cannot be greater than the direct cost minus the entrance cost and the exiting cost, we fix that difference as the upper bound. If the price exceeds that value, the follower is discouraged from using the hub backbone network. Finally, we generate a vector of prices by generating values between the lower and upper bounds. Notice that, since we are using only hubs such that $i > j$, there is no need to use the full matrix; only the lower triangular matrix is required. This is illustrated by the values highlighted in Fig. \ref{fig:CodeP}.

\subsubsection{Solving the Follower's Problem}

We are employing a nested framework in the proposed Co-EA, meaning that the follower's problem must be optimally solved for each leader's decision (tree-of-hubs and prices). This framework ensures that we obtain feasible solutions for the BTHLPwP.

Since the procedure for solving the follower's problem will be repeatedly called in the Co-EA algorithm, it is crucial to implement an efficient method for achieving this. Recall that the follower's problem is parameterized by the leader's variables $y$, $z$, and $\pi$. Once these variables are fixed, the follower's problem results in a series of shortest-path problems (as many as commodities), which could be efficiently solved by the well-known Dijkstra's algorithm. With this strategy, the follower's problem is solved and variables $x$, $x^1$, $x^2$, and $q$ are obtained.

\subsection{Evaluating the Individuals' Fitness}

Once an initial population of tree-structured networks ($y$, $z$) and a population of pricing vectors ($\pi$) are constructed, a matching strategy must be defined. As explained above, the creation of a pricing vector is associated with a specific tree. Hence, the matching criterion is straightforward, and the tuple is $(y, z, \pi)$.

The resulting follower's problem for $(y, z, \pi)$ is optimally solved to obtain ($x$, $x^1$, $x^2$, $q$). Then, the leader's objective function $F(z, \pi, x, x^1)$, defined in expression \eqref{l:0}, can be evaluated. In other words, the fitness of individual $(y, z, \pi)$ is measured.

\subsection{Evolution of the Tree-Structured Networks' Population}

A partial leader's solution, consisting of the tree-structured network, evolves through a standard EA, which has demonstrated effectiveness in solving similar problems \citep{yang2020nature}. In this process, individuals with better fitness will have more opportunities to undergo the evolutionary operators to create new individuals. A description of this process is provided below:

\begin{itemize}
    \item \textbf{Initial population}: Individuals are generated as described in Section \ref{sec:constructing}.
    \item \textbf{Selection}: Based on the individuals' fitness, pairs of them are chosen in an elitist and pseudo-random manner.
    \item \textbf{Crossover}: We implement a two-point crossover mechanism to generate offspring. In this mechanism, two parents are selected. Then, two positions in the solution's vector of one parent are randomly identified, and the offspring inherit characteristics from both parents.
    \item \textbf{Repairing Phase}: The crossover operator may lead to infeasible individuals, which are repaired in this phase. Infeasibility is caused by activating a different number of the $p$ required hubs. If the number of activated hubs is less than $p$, we randomly activate a new hub. The process is repeated until exactly $p$ hubs are activated. On the other hand, if the number of activated hubs exceeds $p$, we randomly close a hub. Similarly, the process is repeated until a feasible individual is obtained.
    \item \textbf{Mutation}: To add diversity into the evolutionary process, new trees are considered in the mutation operator. Since the number of different trees with $n$ nodes is $n^{n-2}$, we use a new configuration of arcs connecting the tree.
\end{itemize}

\begin{figure}[H]
    \centering
    \includegraphics[width = \textwidth]{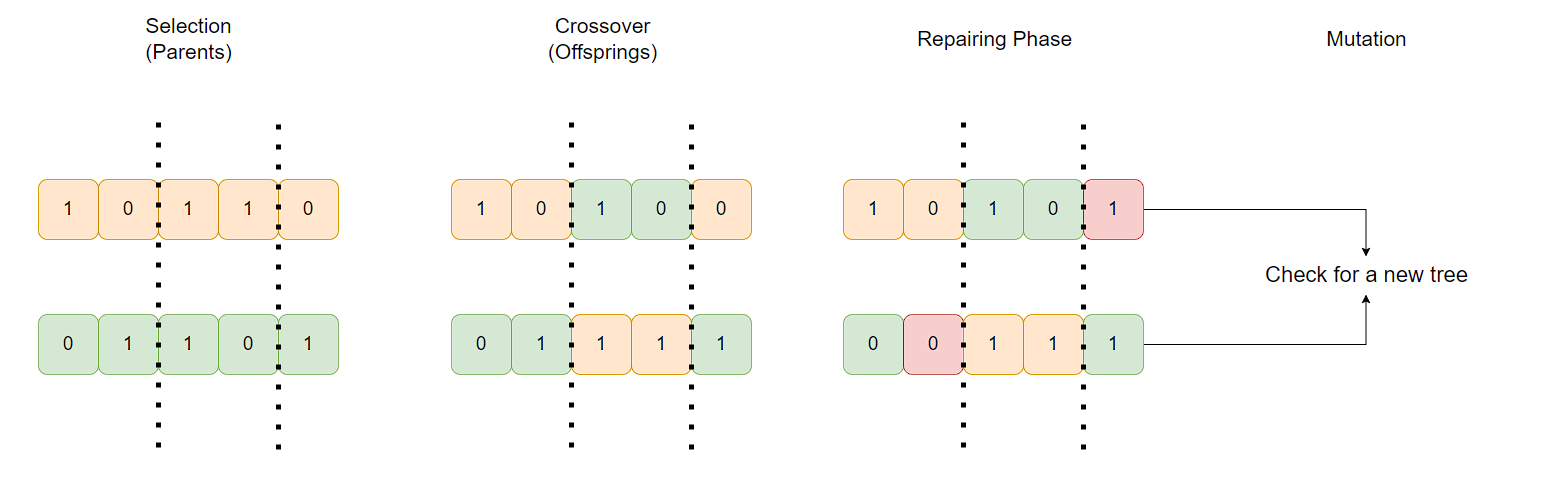}
    \caption{Illustration of the evolutionary operators applied to the tree-structured networks with $p=3$.}
    \label{fig:EAH}
\end{figure}

An illustration of the implemented evolutionary operators for the tree-structured networks is presented in Figure \ref{fig:EAH}.

\subsection{Evolution of the Prices' Population}

Similar to the above, we follow an elitist approach, meaning individuals with better fitness will have better opportunities to survive and crossover. Offspring will inherit characteristics from both parents, and potentially, they may undergo mutation to introduce diversity.

The details of this process are described below:
\begin{itemize}
\item \textbf{Initial population}: Individuals (pricing vectors) are generated in a random manner but within clearly defined bounds, following the procedure described in Section \ref{sec:constructing}.
\item \textbf{Selection}: Fitness is measured according to the leader's objective function. Then, a fixed number of individuals is selected in an elitist manner.
\item \textbf{Crossover}: Pairs of individuals are matched to perform the crossover. Pairs are created by selecting one individual from the elitist set and another from the initial population. A pair consists of two parents and will generate two new individuals (offspring). The crossover is carried out by separating the price vector into quadrants and recombining them to create the new offspring.
\item \textbf{Mutation}: In this phase, the generated offspring are considered. We separate the offspring into quadrants and identify the quadrant with the smallest size. Then, we increase the values in this quadrant by adding a small value.
\end{itemize}

\begin{figure}[H]
    \centering
    \includegraphics[width = \textwidth]{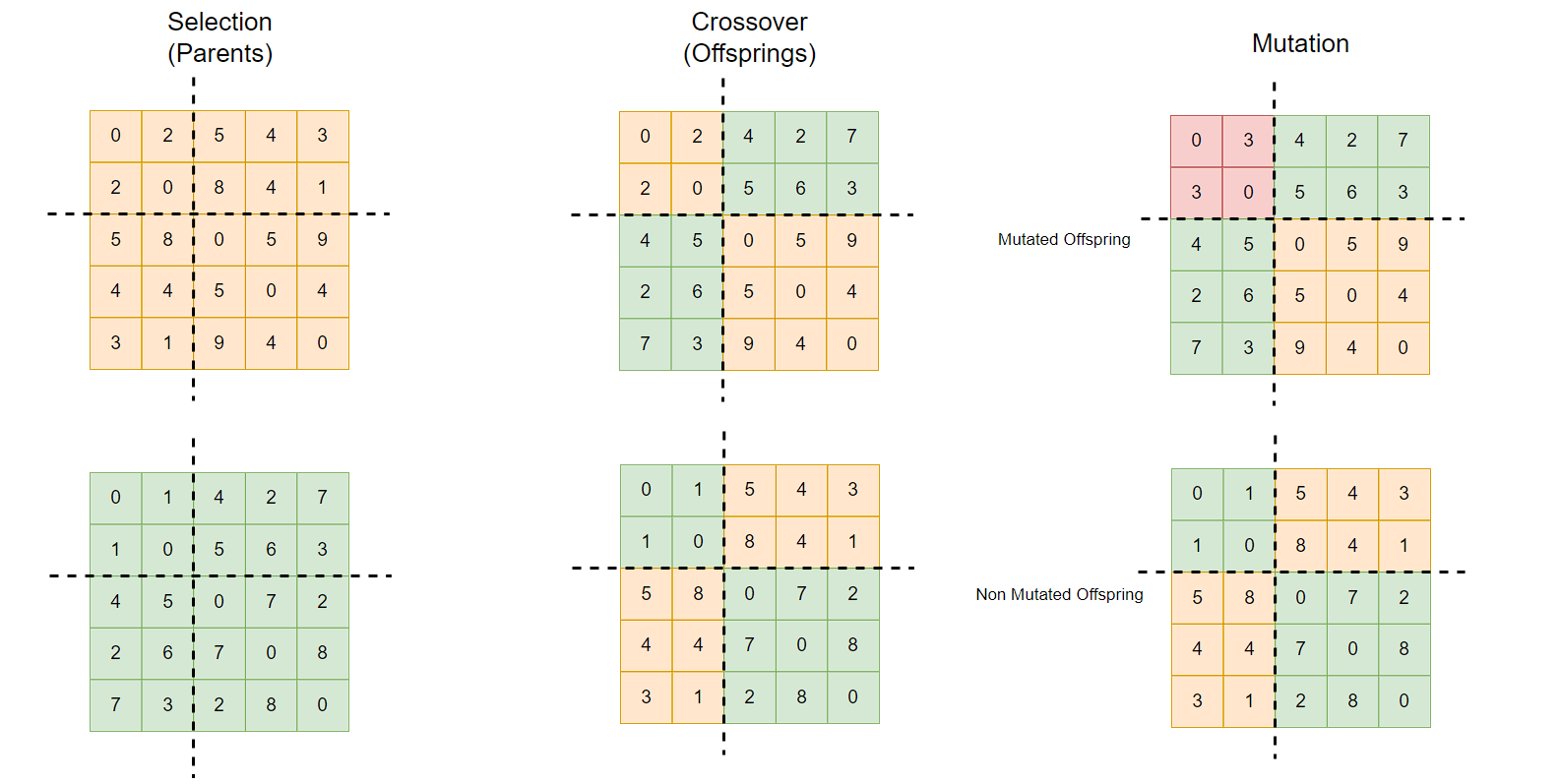}
    \caption{Illustration of the evolutionary operators applied to the pricing vectors with $p=3$.}
    \label{fig:EAP}
\end{figure}

The evolving process of the pricing vectors is depicted in Figure \ref{fig:EAP}.

\subsection{Co-Evolution Between Populations} 
As already mentioned, cooperative co-evolution involves leveraging characteristics of one population to enhance the other. Under this framework, once the population of trees has been updated, it is used to conduct the process again to find bounds on the pricing vectors. Then, the population of the new pricing vectors evolves, and they are used to update the population of tree-structured networks.

Specifically, the new pricing vectors are jointly considered with the current population of trees, and the follower's problem is solved again. By doing this, the hubs with higher frequency (the ones appearing most frequently in solutions) are identified, and that information is used to create new trees. It is evident that there is an information exchange from the pricing population to the population of trees. Furthermore, new trees will be used to compute new bounds on the pricing, leading to a new pricing population. Hence, an information exchange also exists from the population of trees to the pricing population; that is, the co-evolutionary process is performed in two ways. This process is repeated until a stopping criterion is met.

The first stopping criterion of the Co-EA involves reaching, for a second consecutive time, a predetermined number of iterations without reporting an improvement in the leader's objective function. Upon reaching this predetermined number of iterations for the first time, a diversification component is added. Specifically, the current population of trees is replaced by a new population of trees, greedily generated. The algorithm continues its process aiming to improve the solutions until a maximum number of generations is reached. When this occurs, the algorithm ends and reports the individual associated with the best leader's fitness.

\subsection{Parameters in the Co-EA} \label{sec:parameters}

In the proposed algorithm, some parameters must be identified and calibrated to enhance its performance. Specifically, these parameters include the size of the population of trees, the size of the pricing population, the crossover and mutation probabilities, the number of iterations without improvement, and the maximum number of generations. In order to enhance the performance of the Co-EA, these parameters will be calibrated.

\section{Computational Experiments} \label{sec:exp}

{ In this section report the results of our compuational experience. We divided our study into three main parts. First, we implement the single-level reformulation described in Section \ref{sec:SLR} using Gurobi 12.0 to solve small-size instances of the BTHLPwP. Second, we apply the proposed Co-EA to solve all the tested instances. Third, we use the best solution obtained by the Co-EA as a warm start in Gurobi to solve the small-size instances.}

For our experiments, we adapted the well-known AP, CAB, and TR datasets for hub location problems~\citep{wandelt2022toward}. The sizes of the instances are shown in Table \ref{tabla: tabla0}. We selected instances with $10$, $15$, $20$, and $25$ potential sites and customers. For instances with $10$ potential sites and customers, $3$ and $5$ hubs must be located. For the remaining instances ($15$, $20$, and $25$ potential sites), we locate $3$, $5$, and $7$ hubs, respectively. Each instance is created using $5$ different construction costs to assess the impact of varying values. With these configurations, we have $55$ instances for each of the three datasets (AP, CAB, and TR), resulting in a total of $165$ instances.

\begin{table}[H]
\begin{center}
\begin{tabular}{ccccc|}
\hline
\multicolumn{1}{|c}{Instances}            & Potential Sites      & Customers            & $p$       & $\#$ Instances \\ \hline
\multicolumn{1}{|c}{\multirow{4}{*}{AP}}  & 10                   & 10                    & 3, 5      & 10             \\
\multicolumn{1}{|c}{}                     & 15                   & 15                    & 3, 5, 7 & 15             \\
\multicolumn{1}{|c}{}                     & 20                   & 20                    & 3, 5, 7 & 15             \\
\multicolumn{1}{|c}{}                     & 25                   & 25                    & 3, 5, 7 & 15             \\ \hline
\multicolumn{1}{|c}{\multirow{4}{*}{CAB}}  & 10                   & 10                    & 3, 5      & 10             \\
\multicolumn{1}{|c}{}                     & 15                   & 15                    & 3, 5, 7 & 15             \\
\multicolumn{1}{|c}{}                     & 20                   & 20                    & 3, 5, 7 & 15             \\
\multicolumn{1}{|c}{}                     & 25                   & 25                    & 3, 5, 7 & 15             \\ \hline
\multicolumn{1}{|c}{\multirow{4}{*}{TR}} & 10                   & 10                    & 3, 5      & 10             \\
\multicolumn{1}{|c}{}                     & 15                   & 15                    & 3, 5, 7 & 15             \\
\multicolumn{1}{|c}{}                     & 20                   & 20                    & 3, 5, 7 & 15             \\
\multicolumn{1}{|c}{}                     & 25                   & 25                    & 3, 5, 7 & 15             \\ \hline
\multicolumn{1}{l}{}                      & \multicolumn{1}{l}{} & \multicolumn{1}{l|}{} & Total     & 165            \\ \cline{4-5} 
\end{tabular}
\caption{Summary of the tested instances.}
\label{tabla: tabla0}
\end{center}
\end{table}

All the computational experiments were carried out using an Intel(R) Core i5-10400F processor with 16 GB of RAM under Windows 10 Professional.

\subsection{{Results Obtained from the Single-Level Reformulation}} \label{sec:resultsSLR}

{Although the single-level reformulation may result in a mixed-integer linear problem, the use of big-M constants and the need to incorporate constraints to avoid subtours complicate obtaining solution with the Branch-and-Bound approach, as mentioned at the end of Section \ref{sec:SLR}. Nevertheless, we implemented it in Gurobi 12.0 to solve the smallest-size instances of each subset. Specifically, a subset of the AP, CAB, and TR datasets is considered, including 10 potential sites and 10 customers with $p=3$ and $p=5$, resulting in a total of 30 instances.

We set maximum time limit of 6 hours in Gurobi, but none of the instances was solved to optimality within this time. Furthermore, the smallest optimality gap for the best feasible solution reported by Gurobi is around $40\%$. The results are included in the left section of Table \ref{tab:comparison}.}

\subsection{Results Obtained Using the Co-EA} \label{sec:resultsCoEA}

Before using the Co-EA to solve the BTHLPwP, the calibration of the involved parameters needs to be discussed. Recall that the parameters used in the Co-EA are summarized in Section \ref{sec:parameters}.
 Preliminary tests were conducted to determine the parameters for our computational experiments. However, it is important to note that \cite{sipper2018investigating} presents results from large-scale experiments focusing on the robustness of parameter settings for EAs, concluding that, in general, parameter settings are not critical for sufficiently long runs.

Based on the preliminary experiments, we set the size of the initial population of trees and prices generated in the proposed Co-EA to 50. Since both populations evolve through an EA, evolutionary operators are applied. The probability of a pair of solutions entering the crossover is set to $0.8$ for both populations. The result of the crossover then enters the mutation phase with a probability of $0.2$ in each population. These values fall within the range recommended by \cite{sipper2018investigating} as appropriate. The algorithm terminates after $40$ iterations or when $15$ consecutive iterations without improvement of the incumbent solution are observed.

The Co-EA is implemented in Python 3.9. Due to its inherent randomness, 10 runs were performed for each instance. The tag $n$-$p$-\texttt{l} indicates the number of potential sites ($n$), which coincides with the number of users in all instances, the number of hubs to locate ($p$), and an instance identifier (\texttt{l}). We created five different instances of each size, labeled \texttt{A}, \texttt{B}, \texttt{C}, \texttt{D}, and \texttt{E}. We denote by $F_{\text{best}}$ the best leader's profit obtained in the 10 runs. The average and the worst profits are denoted as $F_{\text{avg}}$ and $F_{\text{worst}}$, respectively. Additionally, \texttt{CPUTime} reports the average required computational time (in seconds) for the 10 runs. The detailed results of our experiments are shown in Tables \ref{tab:CABDataset}, \ref{tab:TRDataset}, and \ref{tab:APDataset} in \ref{app:ce}.

To improve the readability of the obtained results, Table \ref{tab: SumDeviation} provides a compact summary showing the different deviations between the leader's profits. This table is organized into three sections, each corresponding to the different types of instances used in the experiments (CAB, AP, and TR). For each set of instances, we separated the results by the number of potential sites ($n$) and the number of hubs to locate ($p$). In the \textbf{Dev. Best vs. Average} column, the deviation between the best profit found and the average profit over the 10 runs is calculated. Likewise, we computed the deviation between the best and the worst profits found in the 10 runs, which is represented by \textbf{Dev. Best vs. Worst}. For each of these deviations, the minimum, average, and maximum values obtained across the five instances (A, B, C, D, and E) are reported. In the \textbf{Avg. CPU Time} column, we report the average running times for solving each of the instances.

\begin{table}[h]
\centering
\begin{adjustbox}{max width=0.8\textwidth}
\begin{tabular}{ccc|ccc|ccc|c}
\multirow{2}{*}{Instance} & \multirow{2}{*}{$n$}  & \multirow{2}{*}{$p$} & \multicolumn{3}{c|}{\textbf{Dev. Best vs. Average}} & \multicolumn{3}{c|}{\textbf{Dev. Best vs. Worst}} & \multirow{2}{*}{\textbf{Avg. CPU Time} } \\
                           &                     &                    & Min         & Average      & Max        & Min       & Average     & Max      &                   \\ \hline
 \multirow{11}{*}{CAB}     & \multirow{2}{*}{10} & 3                  & 0.13        & 0.18         & 0.21       & 0.32      & 0.39        & 0.53     & 611.12                       \\
                           &                     & 5                  & 0.07        & 0.16         & 0.25       & 0.21      & 0.57        & 0.38     & 714.14                        \\ \cline{2-10} 
                           & \multirow{3}{*}{15} & 3                  & 0.26        & 0.64         & 1.60       & 0.52      & 1.36        & 3.44     & 1384.05                       \\
                           &                     & 5                  & 0.54        & 1.05         & 1.53       & 1.28      & 2.15        & 4.20     & 1673.41                       \\
                           &                     & 7                  & 0.46        & 0.91         & 1.66       & 1.06      & 1.90        & 3.04     & 2089.76                       \\ \cline{2-10} 
                           & \multirow{3}{*}{20} & 3                  & 0.86        & 3.15         & 6.60       & 1.93      & 6.50        & 11.28    & 2039.86                      \\
                           &                     & 5                  & 0.63        & 2.37         & 6.61       & 1.61      & 4.40        & 12.15    & 2392.64                       \\
                           &                     & 7                  & 0.71        & 2.23         & 3.58       & 1.63      & 4.52        & 6.97     & 2730.39                       \\ \cline{2-10} 
                           & \multirow{3}{*}{25} & 3                  & 2.71        & 4.57         & 7.12       & 5.72      & 9.33        & 13.09    & 3560.84                       \\
                           &                     & 5                  & 3.86        & 5.27         & 8.21       & 8.04      & 10.69       & 13.70    & 4530.15                       \\
                           &                     & 7                  & 1.81        & 5.26         & 7.72       & 4.32      & 10.89       & 14.02    & 4361.13                       \\ \hline
\multirow{11}{*}{AP}       & \multirow{2}{*}{10} & 3                  & 0.17        & 0.21         & 0.26       & 0.32      & 0.47        & 0.54     & 595.57                       \\
                           &                     & 5                  & 0.12        & 0.17         & 0.25       & 0.28      & 0.34        & 0.39     & 907.97                       \\ \cline{2-10} 
                           & \multirow{3}{*}{15} & 3                  & 0.74        & 1.22         & 1.52       & 1.68      & 2.76        & 3.53     & 1220.54                        \\
                           &                     & 5                  & 0.56        & 1.70         & 2.66       & 1.83      & 3.46        & 4.77     & 1312.93                      \\
                           &                     & 7                  & 0.67        & 1.45         & 3.05       & 1.35      & 2.90        & 4.83     & 1577.48                        \\ \cline{2-10} 
                           & \multirow{3}{*}{20} & 3                  & 0.44        & 1.54         & 2.38       & 0.79      & 3.11        & 4.60     & 1824.92                    \\
                           &                     & 5                  & 0.55        & 2.86         & 4.72       & 1.22      & 5.24        & 8.65     & 1881.68                       \\
                           &                     & 7                  & 1.60        & 2.23         & 3.39       & 2.85      & 3.82        & 4.94     & 1933.42                       \\ \cline{2-10} 
                           & \multirow{3}{*}{25} & 3                  & 2.82        & 6.62         & 9.58       & 5.68      & 13.49       & 18.83    & 2488.75                      \\
                           &                     & 5                  & 2.18        & 6.51         & 12.30      & 5.98      & 11.50       & 18.46    & 2714.34                       \\
                           &                     & 7                  & 1.30        & 5.35         & 9.02       & 2.81      & 10.85       & 18.96    & 2586.29                       \\ \hline
\multirow{11}{*}{TR}       & \multirow{2}{*}{10} & 3                  & 0.10        & 0.16         & 0.32       & 0.19      & 0.38        & 0.95     & 686.01                       \\
                           &                     & 5                  & 0.03        & 0.60         & 0.97       & 0.14      & 0.44        & 0.88     & 755.87                        \\ \cline{2-10} 
                           & \multirow{3}{*}{15} & 3                  & 0.32        & 1.12         & 2.12       & 1.19      & 2.02        & 3.23     & 1721.21                       \\
                           &                     & 5                  & 1.28        & 1.90         & 2.39       & 2.06      & 3.03        & 4.06     & 2034.02                        \\
                           &                     & 7                  & 0.95        & 1.52         & 2.14       & 1.57      & 2.45        & 3.43     & 2124.59                      \\ \cline{2-10} 
                           & \multirow{3}{*}{20} & 3                  & 1.26        & 4.15         & 6.43       & 3.58      & 8.73        & 12.15    & 2414.11                       \\
                           &                     & 5                  & 1.38        & 2.61         & 4.13       & 3.05      & 5.29        & 9.02     & 2827.01                       \\
                           &                     & 7                  & 1.83        & 3.53         & 4.99       & 3.70      & 5.96        & 6.81     & 3612.12                       \\ \cline{2-10} 
                           & \multirow{3}{*}{25} & 3                  & 1.43        & 7.86         & 14.64      & 2.68      & 17.70       & 27.97    & 5147.35                       \\
                           &                     & 5                  & 1.22        & 6.66         & 10.49      & 2.48      & 12.20       & 18.09    & 8404.67                       \\
                           &                     & 7                  & 4.19        & 5.07         & 5.89       & 7.94      & 9.69        & 11.58    & 6264.35                       \\ \cline{2-10} 
\end{tabular}
\end{adjustbox}
\caption{Summary of the deviations computed.}
\label{tab: SumDeviation} 
\end{table}

To begin the discussion of the obtained results, it is noteworthy that for all the instances, the convergence of the Co-EA with respect to the network design was faster than that of the pricing decisions. This finding is supported by the reported deviations. {We would like to emphasize that for the smaller instances (10 potential sites), the deviation was below 1\% in all cases except for the TR 10-5 instances.} However, as the size of the instances increased, the combinatorial nature of the tree-shaped network design problem led to an increase in the reported deviations. This is evident in the column associated with the maximum deviation. {This can be attributed to the continuous nature of the pricing decisions, which is harder to converge.} Additionally, the most difficult-to-solve dataset instances were the TR ones, which reported the maximum deviation among all the instances.

\begin{figure}[h!]
 \begin{center}
\includegraphics[width=0.9\textwidth]{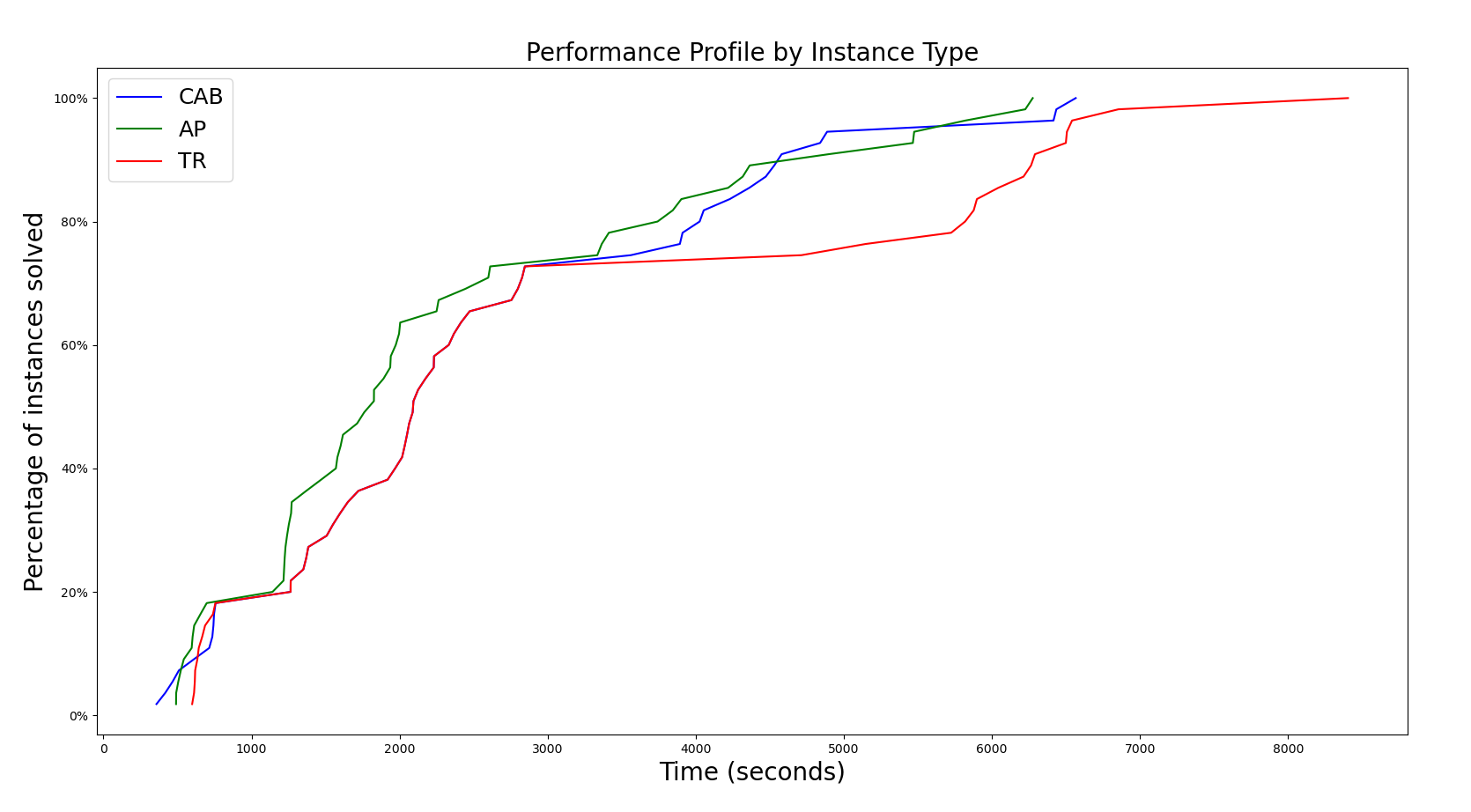}
  \caption{{Performance profile by instance type.}}\label{pp:time}
\end{center}
 \end{figure}

As expected, the computational time increases as the number of potential sites ($n$) and/or the number of hubs to locate ($p$) increases. {However, despite the increase in instance size, the growth in computational time appears to be controlled rather than exponential.} In Figure \ref{pp:time}, we plot the performance profile of the Co-EA, differentiating by the type of instance (CAB, AP, and TR). As observed, instances of type AP seem to be less challenging than those of types CAB and TR, and, furthermore, most instances were solved in less than one hour. {On the other hand, TR instances are the most challenging, with runtimes exceeding two hours. Notably, the cost structure leads to the largest profits among the three types of instances.}

\begin{figure}[h!]
 \begin{center}
\includegraphics[width=\textwidth]{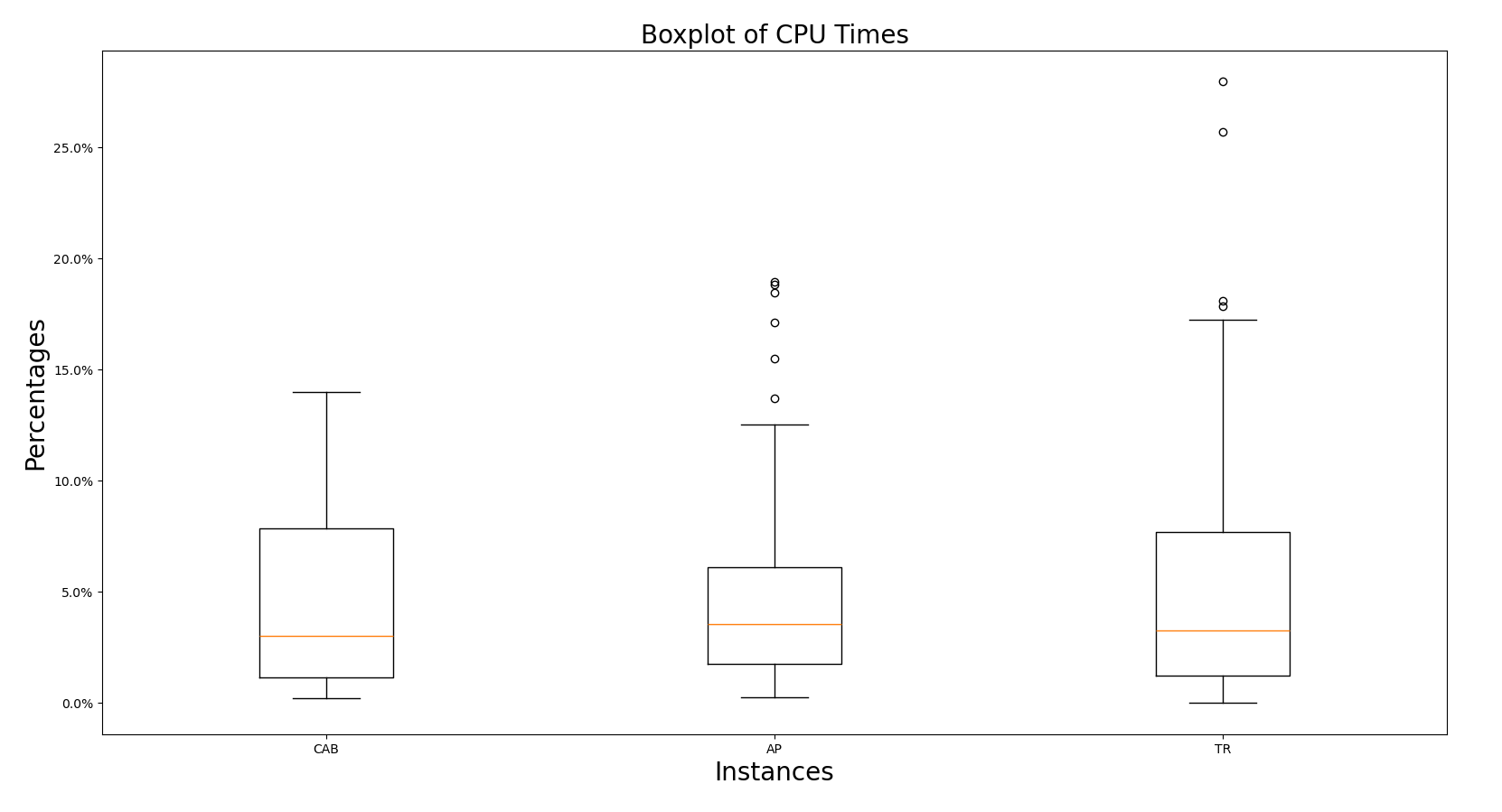}
  \caption{{Boxplots percent deviations reported by the Co-EA in the three datasets of instances.}}\label{boxplot:dev}
\end{center}
 \end{figure}

On the other hand, in Figure \ref{boxplot:dev}, we show the boxplots of the percent deviations of the worst value of each instance with respect to the best obtained profit. As seen in the figure, most deviations fall within the interquartile range of the boxplot. However, there are some atypical cases where the deviations lie outside the box. Although these cases represent large deviations, they occurred only once during the ten runs of the algorithm.

\subsection{{Warm-Starting the Single-Level Reformulation with Solutions Obtained from the Co-EA}}

{Finally, to showcase the benefits of the proposed Co-EA, we use the hub backbone network and prices from the best solution obtained by the Co-EA as a warm start for solving the single-level reformulation of the BTHLPwP using Gurobi 12.0. For this experiment, we use the subset of small-size instances discussed in Section \ref{sec:resultsSLR}. Recall that none of these instances were solved to optimality in 6 hours. In Section \ref{sec:resultsCoEA}, the Co-EA successfully solved these instances and improved all the best-known solutions found by the single-level reformulation. Now, we use the approximate partial solutions (hub backbone network and prices) as input in Gurobi to guide the search process.

It is important to highlight that the hub backbone network ($y, z$) and prices ($\pi$) represent the leader's decision variables. However, in the single-level reformulation, there is no distinction between the leader's and follower's decision variables. As a result, all the primal and dual follower's decision variables, along with the auxiliary variables used for the linearization, are included in the warm-started single-level reformulation. Therefore, to obtain a complete feasible solution, Gurobi must first determine the values of these variables before performing the optimization.

The results of this experiment are presented in Table \ref{tab:comparison}, which is divided into three sections, each corresponding to one solution method: the single-level reformulation in Gurobi, the Co-EA, and the single-level reformulation warm-started with the solution obtained by the Co-EA. Each section includes three columns: the reported leader's profit (in millions), the required computational time (average for the Co-EA) to obtain this value, and a gap column. For the single-level reformulation, the gap represents the optimality gap between the primal and dual bounds computed by Gurobi. For the Co-EA, the gap shows the percentage improvement over the profit found by the single-level reformulation.For the warm-started single-level reformulation, the gap column reports the percentage deviation between the best-known profit found by Gurobi under the warm-start scheme and the best solution obtained by the Co-EA. This column provides information regarding the effectiveness of the Co-EA. For these experiments, a 3-hour time limit was imposed.}

\begin{table}[h]
{ 
\centering
\begin{adjustbox}{max width=0.9\textwidth}
\begin{tabular}{ccc|ccc|ccc|ccc}
\multirow{2}{*}{Instance} & \multirow{2}{*}{}  & & \multicolumn{3}{c|}{\textbf{Single-Level Ref.}} & \multicolumn{3}{c|}{\textbf{Co-EA}} & \multicolumn{3}{c}{\textbf{Warm-Started Single-Level Ref. }} \\
                           &                     &                    & Profit   &  Time (sec)      & Opt. Gap (\%)       & Profit     & Avg. Time (sec)    & Improv. (\%)  & Profit & Time (sec)          & B-K Gap (\%)      \\ \hline
\multirow{10}{*}{CAB}      & 10-3-A &               & 4.79        & $>$21600       & 102.89   & 7.57               & 357.8    & 34.04       & 7.84         & $>$10800       & 3.33    \\
                           &      10-3-B               &                   & 6.90        & $>$21600       & 40.74    & 7.27               & 417.2    & 5.15        & 7.53         & $>$10800       & 3.36    \\ 
                           &      10-3-C               &                   & 4.37        & $>$21600       & 121.80  & 6.00               & 611.12   & 27.15       & 6.29         & $>$10800       & 1.50    \\
                           &     10-3-D                &                   & 6.92        & $>$21600       & 40.60    & 7.78               & 509.19   & 11.04       & 7.80         & $>$10800       & 0.28    \\
                           &     10-3-E                &                   & 4.38        & $>$21600       & 122.09   & 7.49               & 466.66   & 41.59      & 7.70         & $>$10800       & 2.61    \\\cline{2-12} 
                           & 10-5-A &               & 7.97        & $>$21600       & 46.24   & 9.12               & 714.14   & 12.67       & 9.37         & $>$10800       & 2.61    \\
                           &     10-5-B                &                   & 5.75        & $>$21600       & 120.53  & 8.95               & 757.12   & 35.76      & 9.22         & $>$10800       & 2.90   \\
                           &    10-5-C                 &                   & 8.00        & $>$21600       & 45.52    & 9.03               & 746.48   & 11.34       & 9.50         & $>$10800       & 4.98  \\
                           &    10-5-D                 &                   & 6.24        & $>$21600       & 86.86    & 9.53               & 742.49   & 34.51       & 9.72         & $>$10800       & 1.90    \\
                           &    10-5-E                 &                   & 6.83        & $>$21600       & 70.49    & 9.26               & 735.09   & 26.22       & 9.44         & $>$10800       & 1.86   \\\hline
\multirow{10}{*}{AP}      & 10-3-A &                & 3.41        & $>$21600       & 263.01   & 6.39              & 490.07   & 46.63      & 6.66         & $>$10800       & 4.05   \\
                           &   10-3-B                  &                   & 4.74        & $>$21600       & 133.41   & 6.45              & 522.14   & 26.51      & 6.67         & $>$10800       & 3.29    \\ 
                           &    10-3-C                 &                   & 4.72        & $>$21600       & 136.95    & 6.37              & 505.12   & 25.90       & 6.45         & $>$10800       & 1.24   \\
                           &   10-3-D                  &                   & 3.96        & $>$21600       & 216.76    & 6.38              & 490.86  & 37.93        & 6.46         & $>$10800       & 1.23    \\
                           &      10-3-E               &                   & 5.85        & $>$21600       & 91.70     & 6.72              & 595.57   & 12.94       & 7.24         & $>$10800       & 7.18   \\\cline{2-12} 
                           & 10-5-A &               & 6.60        & $>$21600       & 239.54    & 10.93             & 654.23   & 39.61       & 11.25         & $>$10800       & 2.84    \\
                           &     10-5-B                &                   & 8.26        & $>$21600       & 166.37    & 10.80             & 542.47   & 23.51       & 10.91         & $>$10800       & 1.00    \\
                           &    10-5-C                 &                   & 6.91        & $>$21600       & 226.55   & 10.89             & 601.89   & 36.54       & 11.36         & $>$10800       & 4.13    \\
                           &   10-5-D                  &                   & 7.67        & $>$21600       & 186.91   & 10.88             & 612.39   & 29.50       & 11.34         & $>$10800       & 4.05   \\
                           &    10-5-E                 &                   & 8.55        & $>$21600       & 142.55    & 10.86             & 697.12   & 21.27       & 11.32         & $>$10800       & 4.06   \\\hline
\multirow{10}{*}{TR}      & 10-3-A &                & 105.13      & $>$21600       & 138.05    & 152.97            & 738.99   & 31.27       & 153.12        & $>$10800       & 0.09    \\
                           &    10-3-B                 &                   & 90.65       & $>$21600       & 259.48    & 153.11            & 686.01   & 40.79      & 153.73         & $>$10800       & 0.40    \\ 
                           &     10-3-C                &                   & 126.12      & $>$21600       & 88.15     & 152.83            & 598.81   & 17.47      & 152.93         & $>$10800       & 0.06   \\
                           &    10-3-D                 &                   & 117.16      & $>$21600       & 107.89    & 153.08            & 611.76   & 23.46      & 153.19         & $>$10800       & 0.07    \\
                           &    10-3-E                 &                   & 125.54      & $>$21600       & 105.01    & 152.65            & 666.67   & 17.75       & 152.93         & $>$10800       & 0.18   \\\cline{2-12} 
                           & 10-5-A &               & 151.40      & $>$21600       & 206.17    & 258.07            & 634.14   & 41.33       & 258.18         & $>$10800       & 0.04   \\
                           &    10-5-B                 &                   & 171.80      & $>$21600       & 167.10    & 258.14            & 616.71   & 33.44       & 258.73         & $>$10800       & 0.22   \\
                           &    10-5-C                 &                   & 172.30      & $>$21600       & 150.14    & 257.92            & 643.39   & 33.19       & 258.11         & $>$10800       & 0.07    \\
                           &     10-5-D                &                   & 174.51      & $>$21600       & 161.12    & 258.12            & 619.04   & 32.39      & 258.12         & $>$10800       & 0.00   \\
                           &      10-5-E               &                   & 147.77      & $>$21600       & 205.54   & 257.77            & 755.87   & 42.67       & 257.84         & $>$10800       & 0.02    \\\hline
\end{tabular}
\end{adjustbox}
\caption{Showing the benefits of the proposed Co-EA.}
\label{tab:comparison}
}
\end{table}

{From Table \ref{tab:comparison}, it can be concluded that the best profits are obtained by Gurobi under the warm-started single-level reformulation. Nevertheless, it should be noted that these solutions rely on the best partial solution (hubs backbone network and prices) found by the Co-EA. Therefore, the total computational time includes the time required by the Co-EA plus the 3 hours used by Gurobi. Additionally, it is worth mentioning that in all the experiments, the warm-start scheme reached the maximum time limit without finding the optimal solution. 

As mentioned, the aim of this experiment is to validate the effectiveness of the Co-EA, which is successfully demonstrated. In all instances except one, the Co-EA achieves a profit within 5\% of the best-known profit. The only instance with a greater gap is AP 10-3-E, where a gap of 7.18\% is recorded. Nevertheless, it is important to highlight that the best-known solution was obtained using the best solution found by the Co-EA, underscoring the value of the proposed Co-EA. Additionally, it is noteworthy that for the TR instances, a gap of less than or equal to 0.4\% is reported, with instance TR 10-5-D showing that the warm-started scheme did not improve upon the best solution obtained by the Co-EA.
 
It can also be noted from the first section of Table \ref{tab:comparison} that for the CAB instances, the single-level reformulation reported an optimality gap greater than 100\% in four out of the 10 instances. In contrast, the Co-EA obtained a high-quality feasible solution in less than 13 minutes, achieving an improvement of up to 41.5\%. Similarly, for the AP instances, nine out of the 10 instances reported an optimality gap greater than 100\%, while the Co-EA required less than 12 minutes. This behavior remains consistent for the TR instances. These findings highlight the effectiveness and efficiency of the proposed Co-EA in solving this complex problem.}

\subsection{{Sensitivity Analysis}}

{Now, we will explore the significant effect that the number of hubs to be activated has on the profit. To illustrate this analysis, we consider the five CAB instances with 20 potential sites and 3, 5, and 7 hubs to be activated, varying the value of $p$. The profits are plotted in Figure \ref{fig: CABDS20} and detailed in the tables shown in the Appendix.

An interesting observation from our experiments is that, in the BTHLPwP, there is no direct relationship between the number of hubs to be activated and the leader's obtained profit, which contrasts with classical Facility Location problems. Increasing the number of hubs does not necessarily result in a larger profit. While fixed costs increase, the gains from users’ usage may remain similar, thus reducing the overall profit. Therefore, in some instances, we observe that as the value of $p$ increases for the same instance, the profit also increases, while in others, this relationship is not observed.

The blue, red, and green plots in Figure \ref{fig: CABDS20} represent the profit obtained when 3, 5, and 7 hubs are located, respectively. It can be observed that for instances C, D, and E, as the number of hubs increases, the profit also increases. However, in instances A and B, the profit was higher when using 5 hubs instead of 7, contrary to what might have been expected.}

\begin{figure}[H]
    \centering
    \includegraphics[width = \textwidth]{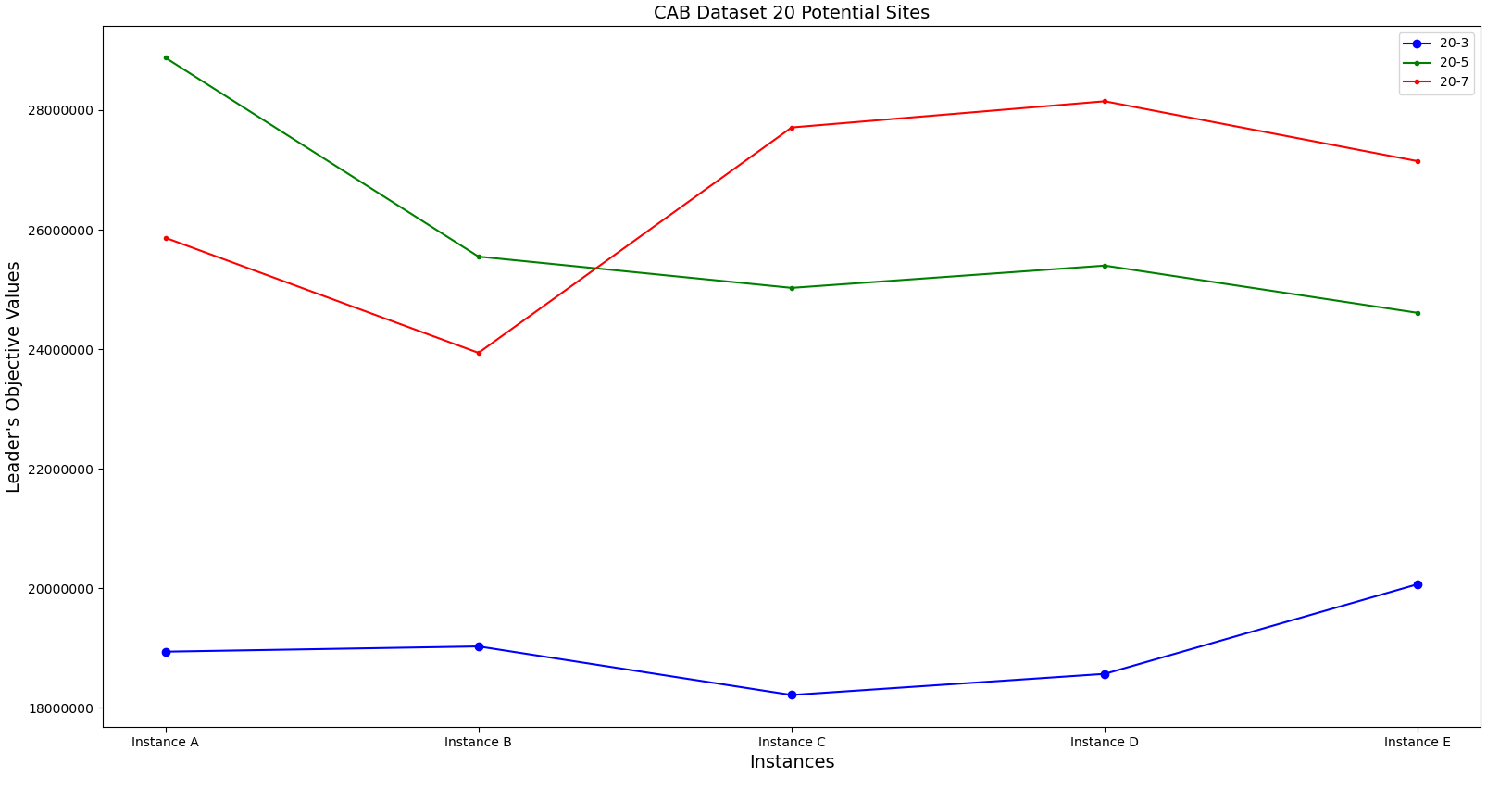}
   \caption{{Leader's profits for the CAB dataset with 20 potential sites, varying the number of located hubs ($p$).}}
    \label{fig: CABDS20}
\end{figure}

{Another finding from the computational experimentation is that, in all instances, the hubs backbone network is used for some commodities. This implies that the proposed algorithm determines a pricing scheme that encourages users to utilize the constructed network. Additionally, some users enter the hubs backbone network and immediately exit without traversing through the arcs. As expected, there are also users who prefer to send their commodities directly, bypassing the constructed network. In these cases, a third-party service is used, but this does not generate revenue for the entity responsible for network design.}

\section{Conclusions and Further Research Directions} \label{sec:concl}

In this study, we introduce the Bilevel Tree of Hubs Location with Pricing, an extension of the Tree of Hubs Location Problem. The development of transportation networks has become essential in modern times, as e-commerce services, airlines, highways, and other industries seek new alternatives for establishing efficient transportation routes. Moreover, understanding the behavior of users utilizing these networks is a crucial factor. {For instance, higher network usage leads to greater profit. Considering this context and the characteristics of the problem, we propose a bilevel model that effectively captures this dynamic. Specifically, we address the problem of designing a tree-structured hub backbone network and setting prices for using its hub links, ensuring that when users route their commodities at minimum cost, the network operator maximizes profit. This general and flexible problem formulation offers valuable insights for practitioners.}

{The problem under study is classified as NP-hard, making it both an interesting and computationally challenging problem. We formulate it as a mixed-integer bilevel optimization problem; however, the lack of specialized solvers for this class of problems necessitates the development of a novel solution methodology. First, we derive an equivalent single-level reformulation by applying the KKT optimality conditions to the follower's problem. This reformulation results in a large-scale mixed-integer linear optimization problem, even for small instances. Although the single-level reformulation can be solved using an optimization solver like Gurobi, it fails to find optimal solutions for small instances within the imposed CPU time limit.
Second, we propose a novel co-evolutionary metaheuristic algorithm that partitions the leader's decisions into two components: constructing the hub backbone network and setting the pricing policy. The algorithm facilitates an innovative exchange of information between these two decision components, specifically tailored for bilevel problems.} It also balances stochasticity and elitism in the generated populations, enhancing solution quality. To evaluate its performance, we adapt and test three different sets of instances from the literature. The proposed algorithm demonstrates consistent performance across various datasets and adapts effectively to different instance characteristics. More broadly, the results highlight the advantages of using co-evolutionary algorithms in settings where a decision-maker's choices interact with those of another decision-maker.
{Third, the solutions obtained by the metaheuristic serve as a valuable warm-start for the exact single-level reformulation, improving the solver's ability to tackle the instances.} It is important to highlight that designing transportation networks is a highly complicated task. So, making long runs of algorithms may lead to better solutions of the considered problem. However, computational time is also an important issue to take into account. {We obtained an efficient reliable solution approach. Thus, we derive a computational experience with the standard datasets in hub location problems, to compare the performance of three procedures based on the developments in this paper. Finally, we use the solution of the co-evolutionary algorithm as warm-start solutions for the single-level reformulation. The main finding is that in three hours, the MIP gaps are reduced, and the quality of the solutions is slightly improved. This implies that the metaheuristic obtains high-quality solutions.}

{As a future research direction, exploring network topologies beyond tree structures could be valuable. While trees typically minimize connection costs, they are vulnerable to disconnection if a node fails, making more reliable topologies an interesting alternative for decision-makers. Another avenue for research is dynamic pricing, which would involve interpreting data to guide its implementation effectively. Additionally, incorporating hub congestion as an objective could further enhance the problem's practical relevance in appropriate applications of this model.}

\section*{Acknowledgements}

The first author acknowledges financial support by  grants PID2020-114594GB-C21 funded by MICIU/AEI/10.13039/501100011033; FEDER+Junta de Andalucía projects C‐EXP‐139‐UGR23, and AT 21\_00032; VII PPIT-US (Ayudas Estancias Breves, Modalidad III.2A); and the IMAG-Maria de Maeztu grant CEX2020-001105-M /AEI /10.13039/501100011033. All the authors acknowledge financial support by grant RED2022-134149-T funded by MICIU/AEI/10.13039/501100011033 (Thematic Network on Location Science and Related Problems).


\begin{thebibliography}{70}
\expandafter\ifx\csname natexlab\endcsname\relax\def\natexlab#1{#1}\fi
\providecommand{\url}[1]{\texttt{#1}}
\providecommand{\href}[2]{#2}
\providecommand{\path}[1]{#1}
\providecommand{\DOIprefix}{doi:}
\providecommand{\ArXivprefix}{arXiv:}
\providecommand{\URLprefix}{URL: }
\providecommand{\Pubmedprefix}{pmid:}
\providecommand{\doi}[1]{\href{http://dx.doi.org/#1}{\path{#1}}}
\providecommand{\Pubmed}[1]{\href{pmid:#1}{\path{#1}}}
\providecommand{\bibinfo}[2]{#2}
\ifx\xfnm\relax \def\xfnm[#1]{\unskip,\space#1}\fi
\bibitem[{Alumur~Alev and Kara(2008)}]{HLP1}
\bibinfo{author}{Alumur~Alev, S.}, \bibinfo{author}{Kara, B.},
  \bibinfo{year}{2008}.
\newblock \bibinfo{title}{Network hub location problems: The state of the art}.
\newblock \bibinfo{journal}{European Journal of Operational Research}
  \bibinfo{volume}{190}, \bibinfo{pages}{1--21}.
\bibitem[{Bagherinejad et~al.(2020)Bagherinejad, Bashiri, Abedpour and
  Soltani}]{bagherinejad2020dynamic}
\bibinfo{author}{Bagherinejad, J.}, \bibinfo{author}{Bashiri, M.},
  \bibinfo{author}{Abedpour, Z.}, \bibinfo{author}{Soltani, E.},
  \bibinfo{year}{2020}.
\newblock \bibinfo{title}{Dynamic single allocation hub location problem
  considering life cycle and reconstruction hubs}.
\newblock \bibinfo{journal}{Research in Production and Operations Management}
  \bibinfo{volume}{11}, \bibinfo{pages}{71--87}.
\bibitem[{Bard(2013)}]{bard2013practical}
\bibinfo{author}{Bard, J.F.}, \bibinfo{year}{2013}.
\newblock \bibinfo{title}{Practical bilevel optimization: algorithms and
  applications}. volume~\bibinfo{volume}{30}.
\newblock \bibinfo{publisher}{Springer Science \& Business Media}.
\bibitem[{Blanco et~al.(2023)Blanco, Fern{\'a}ndez and
  Hinojosa}]{blanco2023hub}
\bibinfo{author}{Blanco, V.}, \bibinfo{author}{Fern{\'a}ndez, E.},
  \bibinfo{author}{Hinojosa, Y.}, \bibinfo{year}{2023}.
\newblock \bibinfo{title}{Hub location with protection under interhub link
  failures}.
\newblock \bibinfo{journal}{INFORMS Journal on Computing} \bibinfo{volume}{35},
  \bibinfo{pages}{966--985}.
\bibitem[{Blanco and Mar{\'i}n(2018)}]{victor}
\bibinfo{author}{Blanco, V.}, \bibinfo{author}{Mar{\'i}n, A.},
  \bibinfo{year}{2018}.
\newblock \bibinfo{title}{Upgrading nodes in tree-shaped hub location}.
\newblock \bibinfo{journal}{Computers and Operations Research}
  \bibinfo{volume}{102}.
\bibitem[{Blanco and Puerto(2022)}]{blanco2022hub}
\bibinfo{author}{Blanco, V.}, \bibinfo{author}{Puerto, J.},
  \bibinfo{year}{2022}.
\newblock \bibinfo{title}{On hub location problems in geographically flexible
  networks}.
\newblock \bibinfo{journal}{International Transactions in Operational Research}
  \bibinfo{volume}{29}, \bibinfo{pages}{2226--2249}.
\bibitem[{Boland et~al.(2004)Boland, Krishnamoorthy, Ernst and
  Ebery}]{boland2004preprocessing}
\bibinfo{author}{Boland, N.}, \bibinfo{author}{Krishnamoorthy, M.},
  \bibinfo{author}{Ernst, A.T.}, \bibinfo{author}{Ebery, J.},
  \bibinfo{year}{2004}.
\newblock \bibinfo{title}{Preprocessing and cutting for multiple allocation hub
  location problems}.
\newblock \bibinfo{journal}{European Journal of Operational Research}
  \bibinfo{volume}{155}, \bibinfo{pages}{638--653}.
\bibitem[{Calvete and Gal{\'e}(2007)}]{calvete2007linear}
\bibinfo{author}{Calvete, H.I.}, \bibinfo{author}{Gal{\'e}, C.},
  \bibinfo{year}{2007}.
\newblock \bibinfo{title}{Linear bilevel multi-follower programming with
  independent followers}.
\newblock \bibinfo{journal}{Journal of Global Optimization}
  \bibinfo{volume}{39}, \bibinfo{pages}{409--417}.
\bibitem[{Camacho-Vallejo et~al.(2023)Camacho-Vallejo, Corpus and
  Villegas}]{camacho2023metaheuristics}
\bibinfo{author}{Camacho-Vallejo, J.F.}, \bibinfo{author}{Corpus, C.},
  \bibinfo{author}{Villegas, J.G.}, \bibinfo{year}{2023}.
\newblock \bibinfo{title}{Metaheuristics for bilevel optimization: A
  comprehensive review}.
\newblock \bibinfo{journal}{Computers \& Operations Research} ,
  \bibinfo{pages}{106410}.
\bibitem[{Camacho-Vallejo and Garcia-Reyes(2019)}]{articleCO2}
\bibinfo{author}{Camacho-Vallejo, J.F.}, \bibinfo{author}{Garcia-Reyes, C.},
  \bibinfo{year}{2019}.
\newblock \bibinfo{title}{Co-evolutionary algorithms to solve hierarchized
  steiner tree problems in telecommunication networks}.
\newblock \bibinfo{journal}{Applied Soft Computing} \bibinfo{volume}{84},
  \bibinfo{pages}{105718}.
\bibitem[{Campbell(1994)}]{campbell1994integer}
\bibinfo{author}{Campbell, J.F.}, \bibinfo{year}{1994}.
\newblock \bibinfo{title}{Integer programming formulations of discrete hub
  location problems}.
\newblock \bibinfo{journal}{European Journal of Operational Research}
  \bibinfo{volume}{72}, \bibinfo{pages}{387--405}.
\bibitem[{Campbell(1996)}]{campbell1996hub}
\bibinfo{author}{Campbell, J.F.}, \bibinfo{year}{1996}.
\newblock \bibinfo{title}{Hub location and the p-hub median problem}.
\newblock \bibinfo{journal}{Operations research} \bibinfo{volume}{44},
  \bibinfo{pages}{923--935}.
\bibitem[{Church and ReVelle(1974)}]{church1974maximal}
\bibinfo{author}{Church, R.}, \bibinfo{author}{ReVelle, C.},
  \bibinfo{year}{1974}.
\newblock \bibinfo{title}{The maximal covering location problem}, in:
  \bibinfo{booktitle}{Papers of the regional science association},
  \bibinfo{organization}{Springer-Verlag}. pp. \bibinfo{pages}{101--118}.
\bibitem[{Contreras et~al.(2011)Contreras, Cordeau and
  Laporte}]{contreras2011dynamic}
\bibinfo{author}{Contreras, I.}, \bibinfo{author}{Cordeau, J.F.},
  \bibinfo{author}{Laporte, G.}, \bibinfo{year}{2011}.
\newblock \bibinfo{title}{The dynamic uncapacitated hub location problem}.
\newblock \bibinfo{journal}{Transportation Science} \bibinfo{volume}{45},
  \bibinfo{pages}{18--32}.
\bibitem[{Contreras et~al.(2009)Contreras, Fern{\'a}ndez and
  Mart{\'i}n}]{contreras2009tight}
\bibinfo{author}{Contreras, I.}, \bibinfo{author}{Fern{\'a}ndez, E.},
  \bibinfo{author}{Mart{\'i}n, A.}, \bibinfo{year}{2009}.
\newblock \bibinfo{title}{Tight bounds from a path based formulation for the
  tree of hub location problem}.
\newblock \bibinfo{journal}{Computers and Operations Research}
  \bibinfo{volume}{36}, \bibinfo{pages}{3117--3127}.
\bibitem[{Contreras et~al.(2010)Contreras, Fern{\'a}ndez and
  Martín}]{contreras2010tree}
\bibinfo{author}{Contreras, I.}, \bibinfo{author}{Fern{\'a}ndez, E.},
  \bibinfo{author}{Martín, A.}, \bibinfo{year}{2010}.
\newblock \bibinfo{title}{The tree of hubs location problem}.
\newblock \bibinfo{journal}{European Journal of Operational Research}
  \bibinfo{volume}{202}, \bibinfo{pages}{390--400}.
\bibitem[{Contreras and O’Kelly(2019)}]{contreras2019hub}
\bibinfo{author}{Contreras, I.}, \bibinfo{author}{O’Kelly, M.},
  \bibinfo{year}{2019}.
\newblock \bibinfo{title}{Hub location problems}, in:
  \bibinfo{booktitle}{Location science}. \bibinfo{publisher}{Springer}, pp.
  \bibinfo{pages}{327--363}.
\bibitem[{Contreras et~al.(2017)Contreras, Tanash and
  Vidyarthi}]{contreras2017exact}
\bibinfo{author}{Contreras, I.}, \bibinfo{author}{Tanash, M.},
  \bibinfo{author}{Vidyarthi, N.}, \bibinfo{year}{2017}.
\newblock \bibinfo{title}{Exact and heuristic approaches for the cycle hub
  location problem}.
\newblock \bibinfo{journal}{Annals of Operations Research}
  \bibinfo{volume}{258}, \bibinfo{pages}{655--677}.
\bibitem[{Cornu{\'e}jols et~al.(1983)Cornu{\'e}jols, Nemhauser and
  Wolsey}]{Cornue}
\bibinfo{author}{Cornu{\'e}jols, G.}, \bibinfo{author}{Nemhauser, G.},
  \bibinfo{author}{Wolsey, L.}, \bibinfo{year}{1983}.
\newblock \bibinfo{title}{The uncapicitated facility location problem}.
\newblock \bibinfo{type}{Technical Report}. Cornell University Operations
  Research and Industrial Engineering.
\bibitem[{Correia et~al.(2018)Correia, Nickel and Saldanha-da
  Gama}]{correia2018stochastic}
\bibinfo{author}{Correia, I.}, \bibinfo{author}{Nickel, S.},
  \bibinfo{author}{Saldanha-da Gama, F.}, \bibinfo{year}{2018}.
\newblock \bibinfo{title}{A stochastic multi-period capacitated multiple
  allocation hub location problem: Formulation and inequalities}.
\newblock \bibinfo{journal}{Omega} \bibinfo{volume}{74},
  \bibinfo{pages}{122--134}.
\bibitem[{{\v{C}}voki{\'c} et~al.(2016){\v{C}}voki{\'c}, Kochetov and
  Plyasunov}]{cvokic1}
\bibinfo{author}{{\v{C}}voki{\'c}, D.D.}, \bibinfo{author}{Kochetov, Y.A.},
  \bibinfo{author}{Plyasunov, A.V.}, \bibinfo{year}{2016}.
\newblock \bibinfo{title}{A leader-follower hub location problem under fixed
  markups}, in: \bibinfo{booktitle}{Discrete Optimization and Operations
  Research: 9th International Conference, DOOR 2016, Vladivostok, Russia,
  September 19-23, 2016, Proceedings 9}, \bibinfo{organization}{Springer}. pp.
  \bibinfo{pages}{350--363}.
\bibitem[{Davis et~al.(2012)Davis, De~Jong, Vose and
  Whitley}]{davis2012evolutionary}
\bibinfo{author}{Davis, L.D.}, \bibinfo{author}{De~Jong, K.},
  \bibinfo{author}{Vose, M.D.}, \bibinfo{author}{Whitley, L.D.},
  \bibinfo{year}{2012}.
\newblock \bibinfo{title}{Evolutionary algorithms}. volume
  \bibinfo{volume}{111}.
\newblock \bibinfo{publisher}{Springer Science \& Business Media}.
\bibitem[{Dempe(2002)}]{dempe2002foundations}
\bibinfo{author}{Dempe, S.}, \bibinfo{year}{2002}.
\newblock \bibinfo{title}{Foundations of bilevel programming}.
\newblock \bibinfo{publisher}{Springer Science \& Business Media}.
\bibitem[{Ebery et~al.(2000)Ebery, Krishnamoorthy, Ernst and
  Boland}]{ebery2000capacitated}
\bibinfo{author}{Ebery, J.}, \bibinfo{author}{Krishnamoorthy, M.},
  \bibinfo{author}{Ernst, A.}, \bibinfo{author}{Boland, N.},
  \bibinfo{year}{2000}.
\newblock \bibinfo{title}{The capacitated multiple allocation hub location
  problem: Formulations and algorithms}.
\newblock \bibinfo{journal}{European journal of operational research}
  \bibinfo{volume}{120}, \bibinfo{pages}{614--631}.
\bibitem[{Erdo{\u{g}}an et~al.(2022)Erdo{\u{g}}an, Battarra and
  Rodr{\'\i}guez-Ch{\'\i}a}]{erdougan2022hub}
\bibinfo{author}{Erdo{\u{g}}an, G.}, \bibinfo{author}{Battarra, M.},
  \bibinfo{author}{Rodr{\'\i}guez-Ch{\'\i}a, A.M.}, \bibinfo{year}{2022}.
\newblock \bibinfo{title}{The hub location and pricing problem}.
\newblock \bibinfo{journal}{European Journal of Operational Research}
  \bibinfo{volume}{301}, \bibinfo{pages}{1035--1047}.
\bibitem[{Esmaeili and Sedehzade(2017)}]{esmaeili}
\bibinfo{author}{Esmaeili, M.}, \bibinfo{author}{Sedehzade, S.},
  \bibinfo{year}{2017}.
\newblock \bibinfo{title}{Designing a hub location and pricing network in a
  competitive environment}.
\newblock \bibinfo{journal}{Journal of Industrial and Management Optimization}
  \bibinfo{volume}{13}, \bibinfo{pages}{1--15}.
\bibitem[{Eydi et~al.(2025)Eydi, Vaziri and Khaleghi}]{eydi2025hierarchical}
\bibinfo{author}{Eydi, A.}, \bibinfo{author}{Vaziri, P.},
  \bibinfo{author}{Khaleghi, A.}, \bibinfo{year}{2025}.
\newblock \bibinfo{title}{Hierarchical p-hub center problem for perishable
  goods}.
\newblock \bibinfo{journal}{Journal of Industrial and Management Optimization}
  \bibinfo{volume}{21}, \bibinfo{pages}{707--730}.
\bibitem[{Farahani et~al.(2013)Farahani, Hekmatfar, Arabani and
  Nikbakhsh}]{farahani2013hub}
\bibinfo{author}{Farahani, R.Z.}, \bibinfo{author}{Hekmatfar, M.},
  \bibinfo{author}{Arabani, A.B.}, \bibinfo{author}{Nikbakhsh, E.},
  \bibinfo{year}{2013}.
\newblock \bibinfo{title}{Hub location problems: A review of models,
  classification, solution techniques, and applications}.
\newblock \bibinfo{journal}{Computers \& Industrial Engineering}
  \bibinfo{volume}{64}, \bibinfo{pages}{1096--1109}.
\bibitem[{da~Gra{\c{c}}a~Costa et~al.(2008)da~Gra{\c{c}}a~Costa, Captivo and
  Cl{\'\i}maco}]{da2008capacitated}
\bibinfo{author}{da~Gra{\c{c}}a~Costa, M.}, \bibinfo{author}{Captivo, M.E.},
  \bibinfo{author}{Cl{\'\i}maco, J.}, \bibinfo{year}{2008}.
\newblock \bibinfo{title}{Capacitated single allocation hub location
  problem—a bi-criteria approach}.
\newblock \bibinfo{journal}{Computers \& operations research}
  \bibinfo{volume}{35}, \bibinfo{pages}{3671--3695}.
\bibitem[{Hakimi(1983)}]{hakimi1983locating}
\bibinfo{author}{Hakimi, S.L.}, \bibinfo{year}{1983}.
\newblock \bibinfo{title}{On locating new facilities in a competitive
  environment}.
\newblock \bibinfo{journal}{European Journal of Operational Research}
  \bibinfo{volume}{12}, \bibinfo{pages}{29--35}.
\bibitem[{Hekmatfar and Pishvaee(2009)}]{hekmatfar2009hub}
\bibinfo{author}{Hekmatfar, M.}, \bibinfo{author}{Pishvaee, M.},
  \bibinfo{year}{2009}.
\newblock \bibinfo{title}{Hub location problem}, in:
  \bibinfo{booktitle}{Facility Location}. \bibinfo{publisher}{Springer}, pp.
  \bibinfo{pages}{243--270}.
\bibitem[{Kalashnikov et~al.(2015)Kalashnikov, Dempe, P{\'e}rez-Vald{\'e}s,
  Kalashnykova and Camacho-Vallejo}]{kalashnikov2015bilevel}
\bibinfo{author}{Kalashnikov, V.V.}, \bibinfo{author}{Dempe, S.},
  \bibinfo{author}{P{\'e}rez-Vald{\'e}s, G.A.}, \bibinfo{author}{Kalashnykova,
  N.I.}, \bibinfo{author}{Camacho-Vallejo, J.F.}, \bibinfo{year}{2015}.
\newblock \bibinfo{title}{Bilevel programming and applications}.
\newblock \bibinfo{journal}{Mathematical Problems in Engineering}
  \bibinfo{volume}{2015}.
\bibitem[{Kay{\i}{\c{s}}o{\u{g}}lu and
  Akg{\"u}n(2021)}]{kayicsouglu2021multiple}
\bibinfo{author}{Kay{\i}{\c{s}}o{\u{g}}lu, B.}, \bibinfo{author}{Akg{\"u}n,
  {\.I}.}, \bibinfo{year}{2021}.
\newblock \bibinfo{title}{Multiple allocation tree of hubs location problem for
  non-complete networks}.
\newblock \bibinfo{journal}{Computers \& Operations Research}
  \bibinfo{volume}{136}, \bibinfo{pages}{105478}.
\bibitem[{Khaleghi and Eydi(2021)}]{khaleghi2021robust}
\bibinfo{author}{Khaleghi, A.}, \bibinfo{author}{Eydi, A.},
  \bibinfo{year}{2021}.
\newblock \bibinfo{title}{Robust sustainable multi-period hub location
  considering uncertain time-dependent demand}.
\newblock \bibinfo{journal}{RAIRO-Operations Research} \bibinfo{volume}{55},
  \bibinfo{pages}{3541--3574}.
\bibitem[{Khaleghi and Eydi(2022)}]{khaleghi2022multi}
\bibinfo{author}{Khaleghi, A.}, \bibinfo{author}{Eydi, A.},
  \bibinfo{year}{2022}.
\newblock \bibinfo{title}{Multi-period hub location problem: a review}.
\newblock \bibinfo{journal}{RAIRO-operations Research} \bibinfo{volume}{56},
  \bibinfo{pages}{2751--2765}.
\bibitem[{Kim and O'Kelly(2009)}]{articlemax}
\bibinfo{author}{Kim, H.}, \bibinfo{author}{O'Kelly, M.}, \bibinfo{year}{2009}.
\newblock \bibinfo{title}{Reliable p-hub location problems in telecommunication
  networks}.
\newblock \bibinfo{journal}{Geographical Analysis} \bibinfo{volume}{41},
  \bibinfo{pages}{283 -- 306}.
\bibitem[{Kleinert et~al.(2020)Kleinert, Labb{\'e}, Plein and
  Schmidt}]{kleinert2020there}
\bibinfo{author}{Kleinert, T.}, \bibinfo{author}{Labb{\'e}, M.},
  \bibinfo{author}{Plein, F.a.}, \bibinfo{author}{Schmidt, M.},
  \bibinfo{year}{2020}.
\newblock \bibinfo{title}{There’s no free lunch: on the hardness of choosing
  a correct big-m in bilevel optimization}.
\newblock \bibinfo{journal}{Operations research} \bibinfo{volume}{68},
  \bibinfo{pages}{1716--1721}.
\bibitem[{Kleinert and Schmidt(2023)}]{kleinert2023there}
\bibinfo{author}{Kleinert, T.}, \bibinfo{author}{Schmidt, M.},
  \bibinfo{year}{2023}.
\newblock \bibinfo{title}{Why there is no need to use a big-m in linear bilevel
  optimization: A computational study of two ready-to-use approaches}.
\newblock \bibinfo{journal}{Computational Management Science}
  \bibinfo{volume}{20}, \bibinfo{pages}{3}.
\bibitem[{Korani and Eydi(2021)}]{korani2021bi}
\bibinfo{author}{Korani, E.}, \bibinfo{author}{Eydi, A.}, \bibinfo{year}{2021}.
\newblock \bibinfo{title}{Bi-level programming model and kkt penalty function
  solution approach for reliable hub location problem}.
\newblock \bibinfo{journal}{Expert systems with applications}
  \bibinfo{volume}{184}, \bibinfo{pages}{115505}.
\bibitem[{Labb{\'e} and Yaman(2008)}]{labbe2008solving}
\bibinfo{author}{Labb{\'e}, M.}, \bibinfo{author}{Yaman, H.},
  \bibinfo{year}{2008}.
\newblock \bibinfo{title}{Solving the hub location problem in a star--star
  network}.
\newblock \bibinfo{journal}{Networks: An International Journal}
  \bibinfo{volume}{51}, \bibinfo{pages}{19--33}.
\bibitem[{Luer-Villagra and Marianov(2013)}]{luer}
\bibinfo{author}{Luer-Villagra, A.}, \bibinfo{author}{Marianov, V.},
  \bibinfo{year}{2013}.
\newblock \bibinfo{title}{A competitive hub location and pricing problem}.
\newblock \bibinfo{journal}{European Journal of Operational Research}
  \bibinfo{volume}{231}.
\bibitem[{Marianov et~al.(1999)Marianov, Serra and
  ReVelle}]{marianov1999location}
\bibinfo{author}{Marianov, V.}, \bibinfo{author}{Serra, D.},
  \bibinfo{author}{ReVelle, C.}, \bibinfo{year}{1999}.
\newblock \bibinfo{title}{Location of hubs in a competitive environment}.
\newblock \bibinfo{journal}{European Journal of Operational Research}
  \bibinfo{volume}{114}, \bibinfo{pages}{363--371}.
\bibitem[{Merakl{\i} and Yaman(2017)}]{merakli2017capacitated}
\bibinfo{author}{Merakl{\i}, M.}, \bibinfo{author}{Yaman, H.},
  \bibinfo{year}{2017}.
\newblock \bibinfo{title}{A capacitated hub location problem under hose demand
  uncertainty}.
\newblock \bibinfo{journal}{Computers \& Operations Research}
  \bibinfo{volume}{88}, \bibinfo{pages}{58--70}.
\bibitem[{Mohammadi et~al.(2019)Mohammadi, Jula and
  Tavakkoli-Moghaddam}]{mohammadi2019reliable}
\bibinfo{author}{Mohammadi, M.}, \bibinfo{author}{Jula, P.},
  \bibinfo{author}{Tavakkoli-Moghaddam, R.}, \bibinfo{year}{2019}.
\newblock \bibinfo{title}{Reliable single-allocation hub location problem with
  disruptions}.
\newblock \bibinfo{journal}{Transportation Research Part E: Logistics and
  Transportation Review} \bibinfo{volume}{123}, \bibinfo{pages}{90--120}.
\bibitem[{Nickel et~al.(2001)Nickel, Sch{\"o}bel and Sonneborn}]{nickel2001hub}
\bibinfo{author}{Nickel, S.}, \bibinfo{author}{Sch{\"o}bel, A.},
  \bibinfo{author}{Sonneborn, T.}, \bibinfo{year}{2001}.
\newblock \bibinfo{title}{Hub location problems in urban traffic networks}, in:
  \bibinfo{booktitle}{Mathematical methods on optimization in transportation
  systems}. \bibinfo{publisher}{Springer}, pp. \bibinfo{pages}{95--107}.
\bibitem[{Num et~al.(2011)Num, Legillon, Liefooghe and Talbi}]{articleCO1}
\bibinfo{author}{Num, T.}, \bibinfo{author}{Legillon, F.},
  \bibinfo{author}{Liefooghe, A.}, \bibinfo{author}{Talbi, E.G.},
  \bibinfo{year}{2011}.
\newblock \bibinfo{title}{Cobra: A coevolutionary meta-heuristic for bi-level
  optimization}.
\newblock \bibinfo{journal}{Studies in Computational Intelligence}
  \bibinfo{volume}{482}.
\bibitem[{O'Kelly(1986)}]{articleOKelly1}
\bibinfo{author}{O'Kelly, M.}, \bibinfo{year}{1986}.
\newblock \bibinfo{title}{The location of interacting hub facilities}.
\newblock \bibinfo{journal}{Transportation Science} \bibinfo{volume}{20},
  \bibinfo{pages}{92--106}.
\bibitem[{O'Kelly(1987)}]{articleOKelly2}
\bibinfo{author}{O'Kelly, M.}, \bibinfo{year}{1987}.
\newblock \bibinfo{title}{A quadratic integer program for the location of
  interacting hub facilities}.
\newblock \bibinfo{journal}{European Journal of Operational Research}
  \bibinfo{volume}{32}, \bibinfo{pages}{393--404}.
\bibitem[{O'Kelly et~al.(2015)O'Kelly, Campbell, de~Camargo and
  de~Miranda~Jr}]{o2015multiple}
\bibinfo{author}{O'Kelly, M.E.}, \bibinfo{author}{Campbell, J.F.},
  \bibinfo{author}{de~Camargo, R.S.}, \bibinfo{author}{de~Miranda~Jr, G.},
  \bibinfo{year}{2015}.
\newblock \bibinfo{title}{Multiple allocation hub location model with fixed arc
  costs}.
\newblock \bibinfo{journal}{Geographical Analysis} \bibinfo{volume}{47},
  \bibinfo{pages}{73--96}.
\bibitem[{O’Kelly et~al.(2015)O’Kelly, Luna, De~Camargo and
  De~Miranda}]{o2015hub}
\bibinfo{author}{O’Kelly, M.E.}, \bibinfo{author}{Luna, H.P.L.},
  \bibinfo{author}{De~Camargo, R.S.}, \bibinfo{author}{De~Miranda, G.},
  \bibinfo{year}{2015}.
\newblock \bibinfo{title}{Hub location problems with price sensitive demands}.
\newblock \bibinfo{journal}{Networks and Spatial Economics}
  \bibinfo{volume}{15}, \bibinfo{pages}{917--945}.
\bibitem[{Paredis(2000)}]{jan}
\bibinfo{author}{Paredis, J.}, \bibinfo{year}{2000}.
\newblock \bibinfo{title}{Coevolutionary algorithms}.
  chapter~\bibinfo{chapter}{0}.
\newblock pp. \bibinfo{pages}{224--238}.
\bibitem[{Peker and Kara(2015)}]{peker2015p}
\bibinfo{author}{Peker, M.}, \bibinfo{author}{Kara, B.Y.},
  \bibinfo{year}{2015}.
\newblock \bibinfo{title}{The p-hub maximal covering problem and extensions for
  gradual decay functions}.
\newblock \bibinfo{journal}{Omega} \bibinfo{volume}{54},
  \bibinfo{pages}{158--172}.
\bibitem[{Pessoa et~al.(2017)Pessoa, Santos and Resende}]{pessoa}
\bibinfo{author}{Pessoa, L.}, \bibinfo{author}{Santos, A.},
  \bibinfo{author}{Resende, M.}, \bibinfo{year}{2017}.
\newblock \bibinfo{title}{A biased random-key genetic algorithm for the tree of
  hubs location problem}.
\newblock \bibinfo{journal}{Optimization Letters} \bibinfo{volume}{11},
  \bibinfo{pages}{1371{\--}1384}.
\bibitem[{Pineda and Morales(2019)}]{pineda2019solving}
\bibinfo{author}{Pineda, S.}, \bibinfo{author}{Morales, J.M.},
  \bibinfo{year}{2019}.
\newblock \bibinfo{title}{Solving linear bilevel problems using big-ms: Not all
  that glitters is gold}.
\newblock \bibinfo{journal}{IEEE Transactions on Power Systems}
  \bibinfo{volume}{34}, \bibinfo{pages}{2469--2471}.
\bibitem[{Popovici et~al.(2012)Popovici, Bucci, Wiegand and
  De~Jong}]{popovici2012coevolutionary}
\bibinfo{author}{Popovici, E.}, \bibinfo{author}{Bucci, A.},
  \bibinfo{author}{Wiegand, R.}, \bibinfo{author}{De~Jong, E.},
  \bibinfo{year}{2012}.
\newblock \bibinfo{title}{Handbook of natural computing}.
\bibitem[{Potter and De~Jong(1994)}]{potter1994cooperative}
\bibinfo{author}{Potter, M.A.}, \bibinfo{author}{De~Jong, K.A.},
  \bibinfo{year}{1994}.
\newblock \bibinfo{title}{A cooperative coevolutionary approach to function
  optimization}, in: \bibinfo{booktitle}{International Conference on Parallel
  Problem Solving from Nature}, \bibinfo{organization}{Springer}. pp.
  \bibinfo{pages}{249--257}.
\bibitem[{Pozo et~al.(2020)Pozo, Puerto and Ch{\'i}a}]{pozo}
\bibinfo{author}{Pozo, M.}, \bibinfo{author}{Puerto, J.},
  \bibinfo{author}{Ch{\'i}a, A.}, \bibinfo{year}{2020}.
\newblock \bibinfo{title}{The ordered median tree of hubs location problem}.
\newblock \bibinfo{journal}{TOP} .
\bibitem[{Rodriguez-Martin and Salazar-Gonzalez(2008)}]{rodriguez2008solving}
\bibinfo{author}{Rodriguez-Martin, I.}, \bibinfo{author}{Salazar-Gonzalez,
  J.J.}, \bibinfo{year}{2008}.
\newblock \bibinfo{title}{Solving a capacitated hub location problem}.
\newblock \bibinfo{journal}{European Journal of Operational Research}
  \bibinfo{volume}{184}, \bibinfo{pages}{468--479}.
\bibitem[{Martins~de S{\'a} et~al.(2013)Martins~de S{\'a}, Camargo and
  Miranda}]{Martins}
\bibinfo{author}{Martins~de S{\'a}, E.}, \bibinfo{author}{Camargo, R.},
  \bibinfo{author}{Miranda, G.}, \bibinfo{year}{2013}.
\newblock \bibinfo{title}{An improved benders decomposition algorithm for the
  tree of hubs location problem}.
\newblock \bibinfo{journal}{European Journal of Operational Research}
  \bibinfo{volume}{226}, \bibinfo{pages}{185{\--}202}.
\bibitem[{Samaniego~Mena and Novoa-Hern\'andez(2018)}]{articleCO3}
\bibinfo{author}{Samaniego~Mena, E.}, \bibinfo{author}{Novoa-Hern\'andez, P.},
  \bibinfo{year}{2018}.
\newblock \bibinfo{title}{A hybrid approach for solving dynamic bi-level
  optimization problems}.
\newblock \bibinfo{journal}{Computaci\'on y Sistemas} \bibinfo{volume}{22}.
\bibitem[{Sasaki and Fukushima(2001)}]{sasaki}
\bibinfo{author}{Sasaki, M.}, \bibinfo{author}{Fukushima, M.},
  \bibinfo{year}{2001}.
\newblock \bibinfo{title}{Stackelberg hub location problem}.
\newblock \bibinfo{journal}{Journal of the Operations Research Society of
  Japan} \bibinfo{volume}{44}.
\bibitem[{Sipper et~al.(2018)Sipper, Fu, Ahuja and
  Moore}]{sipper2018investigating}
\bibinfo{author}{Sipper, M.}, \bibinfo{author}{Fu, W.}, \bibinfo{author}{Ahuja,
  K.}, \bibinfo{author}{Moore, J.H.}, \bibinfo{year}{2018}.
\newblock \bibinfo{title}{Investigating the parameter space of evolutionary
  algorithms}.
\newblock \bibinfo{journal}{BioData mining} \bibinfo{volume}{11},
  \bibinfo{pages}{1--14}.
\bibitem[{Sloss and Gustafson(2020)}]{sloss20202019}
\bibinfo{author}{Sloss, A.N.}, \bibinfo{author}{Gustafson, S.},
  \bibinfo{year}{2020}.
\newblock \bibinfo{title}{2019 evolutionary algorithms review}.
\newblock \bibinfo{journal}{Genetic programming Theory and practice XVII} ,
  \bibinfo{pages}{307--344}.
\bibitem[{Soylu and Katip(2019)}]{soylu2019multiobjective}
\bibinfo{author}{Soylu, B.}, \bibinfo{author}{Katip, H.}, \bibinfo{year}{2019}.
\newblock \bibinfo{title}{A multiobjective hub-airport location problem for an
  airline network design}.
\newblock \bibinfo{journal}{European Journal of Operational Research}
  \bibinfo{volume}{277}, \bibinfo{pages}{412--425}.
\bibitem[{Topcuoglu et~al.(2005)Topcuoglu, Corut, Ermis and
  Yilmaz}]{topcuoglu2005solving}
\bibinfo{author}{Topcuoglu, H.}, \bibinfo{author}{Corut, F.},
  \bibinfo{author}{Ermis, M.}, \bibinfo{author}{Yilmaz, G.},
  \bibinfo{year}{2005}.
\newblock \bibinfo{title}{Solving the uncapacitated hub location problem using
  genetic algorithms}.
\newblock \bibinfo{journal}{Computers \& Operations Research}
  \bibinfo{volume}{32}, \bibinfo{pages}{967--984}.
\bibitem[{Vasconcelos et~al.(2011)Vasconcelos, Nassi and
  Lopes}]{vasconcelos2011uncapacitated}
\bibinfo{author}{Vasconcelos, A.D.}, \bibinfo{author}{Nassi, C.D.},
  \bibinfo{author}{Lopes, L.A.}, \bibinfo{year}{2011}.
\newblock \bibinfo{title}{The uncapacitated hub location problem in networks
  under decentralized management}.
\newblock \bibinfo{journal}{Computers \& Operations Research}
  \bibinfo{volume}{38}, \bibinfo{pages}{1656--1666}.
\bibitem[{Wandelt et~al.(2022)Wandelt, Dai, Zhang and Sun}]{wandelt2022toward}
\bibinfo{author}{Wandelt, S.}, \bibinfo{author}{Dai, W.},
  \bibinfo{author}{Zhang, J.}, \bibinfo{author}{Sun, X.}, \bibinfo{year}{2022}.
\newblock \bibinfo{title}{Toward a reference experimental benchmark for solving
  hub location problems}.
\newblock \bibinfo{journal}{Transportation Science} \bibinfo{volume}{56},
  \bibinfo{pages}{543--564}.
\bibitem[{Weber(1929)}]{weber1929theory}
\bibinfo{author}{Weber, A.}, \bibinfo{year}{1929}.
\newblock \bibinfo{title}{Theory of the Location of Industries}.
\newblock \bibinfo{publisher}{University of Chicago Press}.
\bibitem[{Yang(2020)}]{yang2020nature}
\bibinfo{author}{Yang, X.S.}, \bibinfo{year}{2020}.
\newblock \bibinfo{title}{Nature-inspired optimization algorithms}.
\newblock \bibinfo{publisher}{Academic Press}.
\bibitem[{Yu and Gen(2010)}]{yu2010introduction}
\bibinfo{author}{Yu, X.}, \bibinfo{author}{Gen, M.}, \bibinfo{year}{2010}.
\newblock \bibinfo{title}{Introduction to evolutionary algorithms}.
\newblock \bibinfo{publisher}{Springer Science \& Business Media}.

\end{thebibliography}

\appendix

\section{Detailed computational results.}\label{app:ce}
In this appendix, the detailed results obtained from the computational experimentation described in Section \ref{sec:exp} are presented. Tables \ref{tab:CABDataset}, \ref{tab:TRDataset}, and \ref{tab:APDataset} contain the instance indicator along with the best, average, and worst profits obtained in the 10 runs of the Co-EA. Additionally, the average time (in seconds) required  to solve each run for all the instances is reported.

\begin{table}[htbp]
{ 
    \centering
    \resizebox{7cm}{!}{
    \begin{tabular}{|c|rrrr|}
    \hline
   \multicolumn{1}{|r|}{} & \multicolumn{4}{c|}{CAB dataset} \\\hline
   \multicolumn{1}{|r|}{} & \multicolumn{1}{c|}{$F_{best}$} & \multicolumn{1}{c|}{$F_{avg}$} & \multicolumn{1}{c|}{$F_{worst}$}   & \multicolumn{1}{l|}{Time} \\
    \hline 
        10-3-A & 7579789.86  & 7569673.41  & 7550847.87  & 357.8 \\ 
        10-3-B & 7277978.87  & 7266106.82  & 7250203.92  & 417.2 \\ 
        10-3-C & 6009770.36  & 5997034.23  & 5978195.89  & 611.12 \\ 
        10-3-D & 7785161.46  & 7768830.46  & 7760021.75  & 509.19 \\ 
        10-3-E & 7499243.92  & 7486263.72  & 7473209.28  & 466.66 \\ \hline
        10-5-A & 9128063.70  & 9105145.39  & 9093104.78  & 714.14 \\ 
        10-5-B & 8952557.97  & 8933934.60  & 8920928.81  & 757.12 \\ 
        10-5-C & 9030787.48  & 9024203.19  & 9011723.64  & 746.48 \\ 
        10-5-D & 9538181.74  & 9525119.64  & 9509612.92  & 742.49 \\ 
        10-5-E & 9268103.90  & 9254889.63  & 9245576.18  & 735.09 \\ \hline
        15-3-A & 21858598.70 & 21508774.62 & 21107368.23 & 1383.22 \\ 
        15-3-B & 20409270.85 & 20355753.84 & 20303000.47 & 1349.67 \\ 
        15-3-C & 23579809.45 & 23441841.01 & 23240093.59 & 1264.34 \\ 
        15-3-D & 21853493.40 & 21792899.71 & 21703149.83 & 1263.45 \\ 
        15-3-E & 19887129.23 & 19788833.57 & 19740653.96 & 1368.98 \\ \hline
        15-5-A & 27163151.65 & 26923084.59 & 26739195.77 & 1549.77 \\ 
        15-5-B & 25583274.17 & 25245731.31 & 25026430.12 & 1721.55 \\ 
        15-5-C & 24983068.60 & 24731360.63 & 24606186.56 & 1507.26 \\ 
        15-5-D & 25507336.98 & 25370724.83 & 25180015.33 & 1597.90 \\ 
        15-5-E & 24630390.71 & 24254597.17 & 23594976.22 & 1650.06 \\ \hline
        15-7-A & 28179268.18 & 28049607.81 & 27880561.91 & 2034.02 \\ 
        15-7-B & 26355102.50 & 26228863.71 & 25936792.36 & 1969.28 \\ 
        15-7-C & 27330429.51 & 26979800.08 & 26692984.36 & 2174.41 \\ 
        15-7-D & 27223736.49 & 27041098.97 & 26814914.93 & 2064.19 \\ 
        15-7-E & 26599249.76 & 26158905.27 & 25790075.15 & 2049.82 \\ \hline
        20-3-A & 18941061.90 & 17690575.86 & 16803651.35 & 1918.88 \\ 
        20-3-B & 19028903.10 & 18422274.97 & 17737439.81 & 2093.41 \\ 
        20-3-C & 18216841.10 & 17779239.83 & 17340323.27 & 2016.73 \\ 
        20-3-D & 18569665.78 & 18410876.59 & 18211960.12 & 2124.59 \\ 
        20-3-E & 20070474.05 & 19532927.51 & 18523824.11 & 2088.11 \\ \hline
        20-5-A & 28875907.94 & 26965882.26 & 25366456.06 & 2230.35 \\ 
        20-5-B & 25549850.97 & 25330623.25 & 24931117.33 & 2331.58 \\ 
        20-5-C & 25026549.94 & 24743154.18 & 24576537.48 & 2414.11 \\ 
        20-5-D & 25399327.45 & 25240439.04 & 24989447.34 & 2232.33 \\ 
        20-5-E & 24608231.01 & 23963526.76 & 23622619.22 & 2366.48 \\ \hline
        20-7-A & 25865063.53 & 25681964.78 & 25442226.07 & 2827.01 \\ 
        20-7-B & 23939417.43 & 23353684.66 & 22956599.62 & 2472.59 \\ 
        20-7-C & 27709890.82 & 27204590.57 & 26744055.67 & 2755.12 \\ 
        20-7-D & 28147934.36 & 27141064.81 & 26348915.62 & 2845.75 \\ 
        20-7-E & 27145335.43 & 26435532.37 & 25253539.88 & 2798.04 \\ \hline
        25-3-A & 71417725.92 & 67459201.94 & 8292765.93 & 3560.84 \\
        25-3-B & 57692278.39 & 55433416.67& 7022188.03 & 4026.15 \\ 
        25-3-C & 49847917.08 & 48496765.74 & 7779608.08 & 4472.61 \\ 
        25-3-D & 57147532.96 & 55122114.20 & 6334730.37 & 3892.51 \\ 
        25-3-E & 54113218.44 & 50262722.67 & 7745726.48 & 3910.67 \\ \hline
        25-5-A & 78529235.66 & 75106345.62 & 9401275.56 & 4530.15 \\ 
        25-5-B & 81321193.08 & 78179173.01 & 9633241.15 & 4053.66 \\ 
        25-5-C & 90962091.93 & 86229611.42 & 9989281.04 & 6415.76 \\ 
        25-5-D & 90321112.02 & 82902659.25 & 9465055.76 & 4886.71 \\ 
        25-5-E & 83194268.62 & 79276170.55 & 10605750.29 & 4229.60 \\ \hline
        25-7-A & 113301804.21 & 104572327.86 & 10801125.49 & 4579.74 \\ 
        25-7-B & 101703768.38 & 99866719.39 & 10912205.21 & 6435.04 \\ 
        25-7-C & 103676959.34 & 99590719.29 & 10738221.07 & 4361.13 \\ 
        25-7-D & 94979699.99 & 90101516.17 &  11083820.40 & 6565.95 \\ 
        25-7-E & 90672141.70 & 90672141.70 & 11016164.98 & 4839.47 \\ \hline
    \end{tabular}}
    \caption{{Profits obtained for the instances in the CAB dataset.}}
    \label{tab:CABDataset}
    }
\end{table}

\begin{table}[htbp]
{
  \centering
  \resizebox{7cm}{!}{
    \begin{tabular}{|c|rrrr|}
    \hline
   \multicolumn{1}{|r|}{} & \multicolumn{4}{c|}{AP dataset} \\\hline
   \multicolumn{1}{|r|}{} & \multicolumn{1}{c|}{$F_{best}$} & \multicolumn{1}{c|}{$F_{avg}$} & \multicolumn{1}{c|}{$F_{worst}$}   & \multicolumn{1}{l|}{Time} \\
    \hline
        10-3-A & 63993978.98 & 6377351.95  & 6360679.26  & 490.07 \\ 
        10-3-B & 6457653.87  & 6446400.23  & 6424264.59  & 522.14 \\ 
        10-3-C & 6378733.93  & 6365642.03  & 6344416.90  & 505.12 \\ 
        10-3-D & 6384296.60  & 6370234.16  & 6364083.10  & 490.86 \\ 
        10-3-E & 6721887.04  & 6707941.61  & 6692675.75  & 595.57 \\ \hline
        10-5-A & 10937157.66 & 10919276.30 & 10902555.67 & 654.23 \\ 
        10-5-B & 10803000.37 & 10775633.66 & 10760995.24 & 542.47 \\ 
        10-5-C & 10898367.12 & 10882296.69 & 10863820.33 & 601.89 \\ 
        10-5-D & 10880937.12 & 10861830.23 & 10838433.66 & 612.39 \\ 
        10-5-E & 10865712.17 & 10852149.25 & 10835274.55 & 697.12 \\ \hline
        15-3-A & 6351080.92  & 6284735.02  & 6244670.33  & 1266.93 \\ 
        15-3-B & 6592170.74  & 6543301.30  & 6456235.92  & 1219.63 \\ 
        15-3-C & 6635130.92  & 6534054.91  & 6407356.56  & 1223.56 \\ 
        15-3-D & 6708631.16  & 6620833.40  & 6501080.92  & 1141.64 \\ 
        15-3-E & 6228558.06  & 6136080.93  & 6008615.56  & 1229.12 \\ \hline
        15-5-A & 9513053.78  & 9259864.66  & 9059425.00  & 1271.02 \\ 
        15-5-B & 9607346.01  & 9415392.99  & 9192170.74  & 1369.19 \\ 
        15-5-C & 9407213.56  & 9354691.91  & 9235130.92  & 1239.54 \\ 
        15-5-D & 9001080.92  & 8775342.56  & 8608631.16  & 1215.72 \\ 
        15-5-E & 8906710.56  & 8838149.91  & 8728558.06  & 1251.67 \\ \hline
        15-7-A & 10230480.85 & 9918199.20  & 9735950.68  & 1579.83 \\ 
        15-7-B & 10107082.51 & 10039660.86 & 9970682.24  & 1569.42 \\ 
        15-7-C & 10273566.06 & 10152097.14 & 10008561.92 & 1470.09 \\ 
        15-7-D & 10500735.92 & 10416184.03 & 10240678.66 & 1601.34 \\ 
        15-7-E & 11053992.12 & 10884520.78 & 10725737.06 & 1617.27 \\ \hline
        20-3-A & 7867264.83  & 7832409.27  & 7805026.85  & 1895.80 \\ 
        20-3-B & 7967357.33  & 7777338.28  & 7655615.96  & 1679.26 \\ 
        20-3-C & 7756820.83  & 7671108.77  & 7571584.71  & 1688.51 \\ 
        20-3-D & 7870617.31  & 7725707.47  & 7508731.85  & 1874.54 \\ 
        20-3-E & 7756708.84  & 7607619.27  & 7487371.33  & 1733.72 \\ \hline
        20-5-A & 10554679.75 & 10497128.85 & 10426351.16 & 1826.89 \\ 
        20-5-B & 11072567.33 & 10781954.13 & 10567706.46 & 1973.79 \\ 
        20-5-C & 11228567.86 & 10698372.70 & 10256820.83 & 1712.55 \\ 
        20-5-D & 11870617.30 & 11640631.75 & 11343716.35 & 1940.19 \\ 
        20-5-E & 12660706.13 & 12093779.85 & 11735615.62 & 1760.79 \\ \hline
        20-7-A & 11902204.16 & 11700261.44 & 11562571.33 & 1826.32 \\ 
        20-7-B & 12072811.74 & 11812789.16 & 11567706.96 & 2003.77 \\ 
        20-7-C & 10584710.85 & 10415323.28 & 10206610.83 & 1936.39 \\ 
        20-7-D & 11620660.81 & 11350767.71 & 11206770.82 & 1891.13 \\ 
        20-7-E & 12660706.33 & 12232048.57 & 12035615.62 & 1996.08 \\ \hline
        25-3-A & 8193242.07  & 7520494.38  & 6650508.29  & 3335.03 \\ 
        25-3-B & 8519230.54  & 8013662.09  & 7059117.17  & 3364.90 \\ 
        25-3-C & 8375082.99  & 8139296.11  & 7899274.61  & 3902.70 \\ 
        25-3-D & 8907133.11  & 8053813.06  & 7525834.21  & 4217.92 \\ 
        25-3-E & 8053129.75  & 7523107.37  & 7225558.94  & 3412.89 \\ \hline
        25-5-A & 12324310.79 & 10808685.67 & 10049180.49 & 4364.38 \\ 
        25-5-B & 11585026.07 & 10835073.93 & 10291135.72 & 3742.07 \\ 
        25-5-C & 12201490.74 & 11538564.04 & 10836508.33 & 6275.56 \\ 
        25-5-D & 10991230.61 & 10751217.20 & 10333508.09 & 5820.87 \\ 
        25-5-E & 14394034.40 & 13507558.39 & 12855721.49 & 4897.26 \\ \hline
        25-7-A & 1312985.20  & 12305399.81 & 11473645.55 & 3844.82 \\ 
        25-7-B & 11619260.47 & 11230917.77 & 10896670.98 & 6225.24 \\ 
        25-7-C & 10459151.74 & 9737909.17  & 9025910.52  & 4317.04 \\ 
        25-7-D & 13236416.64 & 13064068.36 & 12864275.04 & 5465.95 \\ 
        25-7-E & 13557701.41 & 12335021.07 & 10986763.07 & 5474.81 \\ \hline
    \end{tabular}}
    \caption{{Profits obtained for the instances in the AP dataset.}}
    \label{tab:APDataset}
    }
\end{table}%

\begin{table}[htbp]
{
  \centering
  \resizebox{7cm}{!}{
    \begin{tabular}{|c|rrrr|}
    \hline
    \multicolumn{1}{|r|}{} & \multicolumn{4}{c|}{TR dataset} \\\hline
    \multicolumn{1}{|r|}{} & \multicolumn{1}{c|}{$F_{best}$} & \multicolumn{1}{c|}{$F_{avg}$} & \multicolumn{1}{c|}{$F_{worst}$}  & \multicolumn{1}{c|}{Time} \\
    \hline
        10-3-A & 152972143.12 & 152952628.10 & 152939697.96 & 738.99 \\ 
        10-3-B & 153112107.12 & 153093344.42 & 153082414.27 & 686.01 \\ 
        10-3-C & 153839712.14 & 153824053.65 & 153800858.62 & 598.81 \\ 
        10-3-D & 153087165.47 & 153066183.42 & 153045622.28 & 611.76 \\ 
        10-3-E & 152659345.36 & 152654390.47 & 152644787.79 & 666.67 \\ \hline
        10-5-A & 258073639.15 & 258053495.30 & 258031829.79 & 634.14 \\ 
        10-5-B & 258143568.34 & 258118526.52 & 258099823.55 & 616.71 \\ 
        10-5-C & 257927147.68 & 257915756.76 & 257904371.63 & 643.39 \\ 
        10-5-D & 258124459.12 & 258115828.86 & 258102090.67 & 619.04 \\ 
        10-5-E & 257771356.85 & 257751216.50 & 257735882.51 & 755.87 \\ \hline
        15-3-A & 215745523.71 & 211105155.25 & 208705168.12 & 1383.22 \\ 
        15-3-B & 219643002.53 & 217907488.45 & 216757533.47 & 1349.67 \\ 
        15-3-C & 216202917.37 & 215517157.24 & 213623945.07 & 1264.34 \\ 
        15-3-D & 227796827.10 & 225272043.95 & 224192839.23 & 1263.45 \\ 
        15-3-E & 214934564.16 & 212268782.27 & 209040103.47 & 1368.98 \\ \hline
        15-5-A & 322023931.69 & 316518096.87 & 308955913.56 & 1549.77 \\ 
        15-5-B & 322636309.39 & 315772816.64 & 311932332.11 & 1721.55 \\ 
        15-5-C & 330848115.52 & 322945767.59 & 320862160.81 & 1507.26 \\
        15-5-D & 313119405.64 & 309123530.31 & 306658142.70 & 1597.90 \\ 
        15-5-E & 326589307.10 & 320019913.00 & 317745516.15 & 1650.06 \\ \hline
        15-7-A & 343731732.11 & 340482497.88 & 338211105.20 & 2034.02 \\ 
        15-7-B & 340586811.16 & 337355395.09 & 335239967.19 & 1969.28 \\ 
        15-7-C & 354207296.51 & 346902334.05 & 343833601.80 & 2174.41 \\ 
        15-7-D & 339060503.96 & 331788941.73 & 329793552.39 & 2064.19 \\ 
        15-7-E & 347603487.36 & 342379423.50 & 335690975.91 & 2049.82 \\ \hline
        20-3-A & 405202386.60 & 385080679.77 & 358448818.81 & 1918.88 \\ 
        20-3-B & 393148017.17 & 379138056.58 & 363658071.90 & 2093.41 \\ 
        20-3-C & 419043452.97 & 392115815.27 & 368140111.33 & 2016.73 \\ 
        20-3-D & 410975966.41 & 392443570.18 & 374462165.91 & 2124.59 \\ 
        20-3-E & 402996417.03 & 397908399.25 & 388554704.37 & 2088.11 \\ \hline
        20-5-A & 556755498.34 & 540288768.85 & 532579872.76 & 2230.35 \\ 
        20-5-B & 826473473.23 & 815048690.18 & 801274333.42 & 2331.58 \\ 
        20-5-C & 747435750.84 & 732671139.60 & 708623262.11 & 2414.11 \\ 
        20-5-D & 681400891.01 & 653240014.49 & 619969255.54 & 2232.33 \\ 
        20-5-E & 821313280.88 & 799966591.77 & 781314348.94 & 2366.48 \\ \hline
        20-7-A & 541543222.88 & 517482340.18 & 504642404.13 & 2827.01 \\ 
        20-7-B & 540962298.55 & 524961547.66 & 507147090.11 & 2472.59 \\ 
        20-7-C & 540561175.07 & 522039750.95 & 505975314.94 & 2755.12 \\ 
        20-7-D & 529729097.80 & 520047079.58 & 510132384.77 & 2845.75 \\ 
        20-7-E & 538863124.04 & 511984510.53 & 503034672.65 & 2798.04 \\ \hline
        25-3-A & 716135412.27 & 678315606.97 & 623412354.15  & 3225.01 \\ 
        25-3-B & 712730939.27 & 642245984.87  & 585399054.48  & 5724.72 \\ 
        25-3-C & 692369198.38 & 682455403.14  & 673838938.04   & 5147.35 \\ 
        25-3-D & 805649943.14 & 687733323.74  & 580300822.35   & 6290.65 \\ 
        25-3-E & 853355757.46 & 759312603.89  & 633852563.63  & 5876.65 \\ \hline
        25-5-A & 868915387.85 & 858324487.69  & 847368537.48   & 6506.37 \\ 
        25-5-B & 953241408.32 & 855983446.64  & 780780385.10  & 8404.67 \\ 
        25-5-C & 920297061.31 & 823800380.96  & 761512725.89  & 5818.88\\ 
        25-5-D & 944811979.69 & 897534011.75  & 836208259.55  & 6042.38 \\ 
        25-5-E & 953881156.46 & 892854689.87  & 842562168.08  &  5898.55\\ \hline
        25-7-A & 972132541.36 & 941665724.05  & 921630089.72  &  6141.64 \\ 
        25-7-B & 966366316.36 & 925863480.38  & 889594035.07  & 6264.35 \\ 
        25-7-C & 962287408.27 & 918149967.52  & 870838300.35  & 6854.24 \\ 
        25-7-D & 998577564.52 & 939761731.41  & 894021573.52  & 6499.15 \\ 
        25-7-E & 972154321.36 & 903574423.47  & 929794509.82  & 6120.66 \\ \hline
    \end{tabular}}
    \caption{{Profits obtained for the instances in the TR dataset.}}
  \label{tab:TRDataset}%
  }
\end{table}%

\end{document}